\documentclass[11pt]{article}

\renewcommand{\baselinestretch}{1}   
\usepackage{psfrag,cite,amsmath,graphicx,epsfig,latexsym,subfigure,color}
\usepackage[scriptsize,bf]{caption}
\usepackage{color}
\usepackage{amssymb}
\usepackage{setspace}
\usepackage{array}
\usepackage{booktabs}

\def\bn{\hfill \\ \smallskip\noindent}

\def\prox{\mbox{prox}}

\newcommand{\beq}{\begin{equation}}
\newcommand{\eeq}{\end{equation}}
\newcommand{\st}{{\rm s.t.}}

\newcommand{\bfOmega}{{\mbox{\boldmath $\Omega$}}}

\newcommand{\bZ}{{\mbox{\boldmath $Z$}}}

\newcommand{\cO}{{\mbox{$\mathcal{O}$}}}

\newcommand{\cK}{{\mbox{$\mathcal{K}$}}}
\newcommand{\cG}{{\mbox{$\mathcal{G}$}}}

\newcommand{\cV}{{\mbox{$\mathcal{V}$}}}

\begin{document}
\def\pn {\par\smallskip\noindent}
\def \bn {\hfill \\ \smallskip\noindent}
\newcommand{\fs}{f_1,\ldots,f_s}
\newcommand{\f}{\vec{f}}
\newcommand{\hf}{\hat{f}}
\newcommand{\hx}{\hat{x}}
\newcommand{\hy}{\hat{y}}
\newcommand{\hz}{\hat{z}}
\newcommand{\hw}{\hat{w}}
\newcommand{\tw}{\tilde{w}}
\newcommand{\hlambda}{\hat{\lambda}}
\newcommand{\hbeta}{\hat{\beta}}
\newcommand{\tG}{\widetilde{G}}
\newcommand{\tg}{\widetilde{g}}
\newcommand{\barhx}{\bar{\hat{x}}}
\newcommand{\vecx}{x_1,\ldots,x_m}
\newcommand{\xoy}{x\rightarrow y}
\newcommand{\barx}{{\bar x}}
\newcommand{\bary}{{\bar y}}
\newcommand{\hrho}{\widehat{\rho}}
\newtheorem{theorem}{Theorem}[section]
\newtheorem{lemma}{Lemma}[section]
\newtheorem{corollary}{Corollary}[section]
\newtheorem{proposition}{Proposition}[section]
\newtheorem{definition}{Definition}[section]
\newtheorem{claim}{Claim}[section]
\newtheorem{remark}{Remark}[section]

\def\smskip{\par\vskip 5 pt}
\def\proof{\bn {\bf Proof.} }
\def\QED{\hfill{\bf Q.E.D.}\smskip}
\def\qed{\quad{\bf q.e.d.}\smskip}
\newcommand{\cE}{\mathcal{E}}
\newcommand{\cM}{\mathcal{M}}
\newcommand{\cN}{\mathcal{N}}
\newcommand{\cJ}{\mathcal{J}}
\newcommand{\cT}{\mathcal{T}}
\newcommand{\bx}{\mathbf{x}}
\newcommand{\bp}{\mathbf{p}}
\newcommand{\bX}{\mathbf{X}}
\newcommand{\bY}{\mathbf{Y}}
\newcommand{\bP}{\mathbf{P}}
\newcommand{\bD}{\mathbf{D}}
\newcommand{\bA}{\mathbf{A}}
\newcommand{\bB}{\mathbf{B}}
\newcommand{\bfM}{\mathbf{M}}
\newcommand{\bL}{\mathbf{L}}
\newcommand{\bz}{\mathbf{z}}
\newcommand{\cF}{\mathcal{F}}
\newcommand{\bF}{\mathbf{F}}
\newcommand{\cR}{\mathcal{R}}
\newcommand{\bzero}{\mathbf{0}}

\newcommand{\blue}{\color{blue}}
\newcommand{\red}{\color{red}}

\title{Decomposing Linearly Constrained Nonconvex Problems by a Proximal Primal Dual Approach: Algorithms, Convergence, and Applications}
\author{{Mingyi Hong}\thanks{M. Hong is with the Department of Industrial and Manufacturing Systems Engineering (IMSE), Iowa State University, Ames, IA 50011, USA. Email: \texttt{mingyi@iastate.edu}. M. Hong is supported in part by NSF under Grant CCF-1526078 and by AFOSR under grant 15RT0767.}}
\date{April 4th, 2016}
\maketitle
\begin{abstract}
In this paper, we propose a new decomposition approach named the proximal primal dual algorithm (Prox-PDA) for smooth nonconvex linearly constrained optimization problems. The proposed approach is primal-dual based, where the primal step minimizes certain approximation of the  augmented Lagrangian of the problem, and the dual step performs an approximate dual ascent. The approximation used in the primal step is able to decompose the variable blocks, making it possible to obtain simple subproblems by leveraging the problem structures. 
Theoretically, we show that whenever the penalty parameter in the augmented Lagrangian is larger than a given threshold, the Prox-PDA converges  to the set of stationary solutions, globally and in a  sublinear manner (i.e., certain measure of stationarity decreases in the rate of $\cO(1/r)$, where $r$ is the iteration counter). Interestingly, when applying a variant of the Prox-PDA to the problem of distributed nonconvex optimization (over a connected undirected graph), the resulting algorithm coincides with the popular EXTRA algorithm [Shi {\it et al} 2014], which is only known to work in {\it convex} cases. Our analysis implies that EXTRA and its variants converge globally sublinearly to stationary solutions of certain nonconvex distributed optimization problem. There are many possible extensions of the Prox-PDA, and we present one particular extension to certain nonconvex distributed matrix factorization problem. 
\end{abstract}

\section{Introduction}\label{sec:introduction}
Consider the following optimization problem
\begin{align}\label{eq:original}
\min_{x\in\mathbb{R}^N}\; f(x),\quad \st \;Ax=b
\end{align}
where $f(x):\mathbb{R}^{N}\to \mathbb{R}$ is a closed and smooth function (possibly nonconvex); $A\in\mathbb{R}^{M\times N}$ is a rank deficient matrix; $b\in\mathbb{R}^{M}$ is a known vector.  In this paper we propose a first-order primal-dual method for solving such nonconvex problem. Below we list a few applications of the above model. 

\subsection{Motivating Examples}
\noindent{\bf Distributed Optimization Over Networks.} Consider a network consists of $N$ agents who collectively optimize the following problem 
\begin{align}\label{eq:global:consensus}
\min_{y\in\mathbb{R}} \; f(y):=\sum_{i=1}^{N} f_i(y),
\end{align}
where $f_i(y): \mathbb{R}\to\mathbb{R}$ is a function local to agent $i$ (here $y$ is assumed to be scalar for ease of presentation). Suppose the agents are connected by a network defined by an {\it undirected} graph $\cG=\{\mathcal{V}, \mathcal{E}\}$, with $|\mathcal{V}|=N$ vertices and $|\mathcal{E}|=E$ edges. Each agent can only communicate with its immediate neighbors, and it is responsible for optimizing one component function $f_i$. This problem has found applications in various domains such as distributed consensus \cite{tsitsiklis84thesis,Xiao:2007:DAC:1222667.1222952}, distributed communication networking \cite{hao_consensus,liao15semi-async},  distributed and parallel machine learning \cite{Forero11,mateos_dlasso, shalev14proximaldual} and distributed signal processing \cite{Schizas08,hao_consensus}; for more applications we refer the readers to a recent survey \cite{Giannakis15}. 

Define the node-edge incidence matrix $A\in\mathbb{R}^{E\times N}$ as following: if $e\in\mathcal{E}$ and it connects vertex $i$ and $j$ with $i>j$, then $A_{ev}=1$ if $v=i$, $A_{ev}=-1$ if $v=j$ and $A_{ev}=0$ otherwise.  Using this definition, the {\it signed graph Laplacian matrix} $L_{-}\in\mathbb{R}^{N\times N}$ is given by 
$$L_{-} := A^T A.$$
Introduce $N$ local variables $x=[x_1,\cdots, x_N]^T$, and suppose the graph $\{\cV, \cE\}$ is connected. Then it is clear that the following formulation is equivalent to the global consensus problem, which is precisely problem \eqref{eq:original}
\begin{align}\label{eq:global:consensus:equiv}
\min_{x\in\mathbb{R}^N} \; f(x):=\sum_{i=1}^{N} f_i(x_i),\quad \st\; Ax = 0.
\end{align}

\noindent{\bf Multi-Block Linearly Constrained Problem.}
Consider the following multi-block linearly constrained problem
\begin{align}\label{eq:2block}
\min_{\{y_i\in\mathbb{R}^N\}}\quad f(y)\quad \st\; \sum_{i=1}^{K}A_iy_i = b
\end{align}
where $A_i\in\mathbb{R}^{M\times N}$; $y=[y_1; y_2;\cdots, y_K]$. 

Define a new variable $x:= [y_1;y_2; \cdots; y_K]\in\mathbb{R}^{NK}$, and a new matrix $C = [A_1; A_2; \cdots, A_K]\in\mathbb{R}^{M\times NK}$, then the above problem can also be cast as a special case of problem \eqref{eq:original}
\begin{align*}
\min_{x\in\mathbb{R}^{NK}}\quad f(x)\quad \st\; C x =b.
\end{align*}
Such problem, convex or nonconvex, has wide applications in practice, such as in distributed optimization and coordination (more specifically the sharing problem) \cite{BoydADMM11,hong14nonconvex_admm}, robust Principal Component Analysis \cite{Zhou10RPCA} and rate maximization problem in downlink broadcast communication channels \cite{jindal05}. 

\subsection{\bf Literature Review.} 
The Augmented Lagrangian (AL) methods, or the methods of multipliers, pioneered by Hestenes \cite{Hestenes} and Powell \cite{powell69}, is a classical algorithm for solving nonlinear nonconvex constrained optimization problems \cite{bertsekas82,Nocedal99}. Many existing packages such as LANCELOT \cite{conn91, LANCELOT} are implemented based on this method. Recently, due to the need to solve very large scale nonlinear optimization problems, the AL and its variants regain their popularity, see recent developments in \cite{Curtis16adaptive,Friedlander05, Fernandez12} and the references therein. Also reference \cite{Burachik2008} have developed an AL based algorithm for nonconvex nonsmooth optimization, where subgradients of the augmented Lagrangian are used in the primal update. When the problem is convex and the constraints are linear, Lan and Monterio \cite{Lan2015augmented} have analyzed the iteration complexity for the AL method. More specifically, they have characterized the total number of Nesterov's optimal iterations \cite{Nesterov04} that are required to reach high quality primal-dual solutions. However, despite the generality of these methods, it appears that the AL methods does not decompose well over the optimization variables. Further, the AL method, at least in its classical forms, is difficult to be implemented in a distributed manner. 

Recently, a method named  Alternating Direction Method of Multipliers (ADMM), which is closely related to the AL method, has gained tremendous popularity in solving large-scale structured optimization problems; see a recent survey by Boyd {\it et al} \cite{BoydADMM11}. The method originates in early 1970s \cite{ADMMGlowinskiMorroco,ADMMGabbayMercier}, and has since been studied extensively \cite{Eckstein89, EcksteinBertsekas1992, Glow84, bertsekas97, HongLuo2012ADMM,DengYin2013J,HeXuYuan2011}. The main strength of this algorithm is that it is capable of decomposing a large and complicated problem into a series of small and simple subproblems, therefore making the overall algorithm scalable and easy to implement. However, unlike the AL method, the ADMM is only known to work for convex problems, despite its good numerical performance in nonconvex problems such as the nonnegative matrix factorization \cite{zhang10ADMM_NMF, sun14}, phase retrieval \cite{wen12}, distributed matrix factorization \cite{zhang14}, distributed clustering \cite{Forero11},tensor decomposition \cite{Liavas14} and so on. Only very recently, researchers have begun to rigorously investigate the convergence behavior of ADMM for nonocnvex problems. Zhang \cite{zhang10ADMM_QP} have analyzed a class of splitting algorithms (which includes the ADMM as a special case) for a very special class of nonconvex quadratic problems. It is  shown that these type of methods converge to the stationary solutions {\it linearly} when certain condition on the dual stepsize is met. Ames and Hong have provided an analysis for ADMM when applied to certain nonconvex quadratically constrained, $\ell_1$ penalized quadratic problem, which arises in high-dimensional discriminant analysis.  Hong, Luo and Razaviyayn have developed a  three-step approach for using the ADMM and its proximal gradient variant for solving certain nonconvex consensus and sharing problems that arise in distributed and parallel computation. The key idea of the analysis is to adopt the augmented Lagrangian as the merit function to guide the progress of the algorithm. Li and Pong \cite{li14nonconvex}, Wang, Yin and Zeng \cite{wang15nonconvexadmm} have also used  similar analysis steps to analyze different forms of nonconvex problems.  However, despite these recent successes, it appears that the aforementioned works still pose very restrictive assumptions on the problem types in order to achieve  convergence. For example it is not clear whether the ADMM can be used for the distributed optimization problem \eqref{eq:global:consensus:equiv} over an arbitrary connected graph, despite the fact that for convex problem such application is popular and efficient \cite{Giannakis15,mateos_dlasso,chang14distributed,Wei13,Schizas08}. 

In this paper, we answer the following research question: Is it possible to develop augmented Lagrangian-like decomposition schemes for the linearly constrained nonconvex problem \eqref{eq:original}, with global convergence rate guarantee. Ideally, the resulting algorithm should be able to decompose the updates of different variable blocks so that each of its steps can be easily implemented in distributed manner and/or in parallel. Further, it is desirable that the resulting algorithm would have global convergence and rate of convergence guarantee. To this end, we study a primal-dual algorithm, where the primal step minimizes certain approximation of the  augmented Lagrangian of problem \eqref{eq:original}, and the dual step performs an approximate dual ascent. The approximation used in the primal step is able to decompose the variables, making it possible to obtain simple subproblems by leveraging the problem structures.  Theoretically, we show that whenever the penalty parameter in the augmented Lagrangian is larger than a given threshold, the Prox-PDA converges  to the set of stationary solutions, globally and in a  sublinear manner (i.e., certain measure of stationarity decreases in the rate of $\cO(1/r)$, where $r$ is the iteration counter). We also analyze various different extensions of the algorithm, and discuss their applications to the distributed nonconvex optimization problem \eqref{eq:global:consensus:equiv}. 

\section{The Proposed Algorithm} 
The proposed algorithm builds upon the classical augmented Lagrangian method (also known as the method of multipliers) \cite{bertsekas82,powell69}.  Let us introduce the augmented Lagrangian for problem \eqref{eq:original} as
\begin{align}
L_{\beta}(x,\mu) = f(x)+\langle \mu, Ax-b\rangle +\frac{\beta}{2}\|Ax-b\|^2\label{eq:augmented}
\end{align}
where $\mu\in\mathbb{R}^{M}$ is the Lagrangian dual variable; $\beta>0$ is a penalty parameter.  Let $B\in\mathbb{R}^{M\times N}$ be some arbitrary matrix.  Then the steps of the proposed proximal primal-dual algorithm is given in the following table:
\begin{center}
	\fbox{
		\begin{minipage}{5in}
			\smallskip
			\centerline{\bf {Algorithm 1. The Proximal Primal Dual Algorithm (Prox-PDA)}}
			\smallskip
			At iteration $0$, initialize $\mu^0 =0$ and $x^0\in\mathbb{R}^N$.
			
			At each iteration $r+1$, update variables  by:
			\begin{subequations}
				\begin{align}
				x^{r+1}& =\arg\min_{x\in\mathbb{R}^n}\; f(x) +\langle\mu^r, Ax - b\rangle+\frac{\beta}{2}\|Ax-b\|^2+\frac{\beta}{2}\|x-x^r\|^2_{B^T B};  \label{eq:x:update}\\
				\mu^{r+1}& = \mu^r +\beta (A x^{r+1}	- b).			\label{eq:mu:update}
				\end{align}
			\end{subequations}
			
		\end{minipage}
	}
\end{center}
In Prox-PDA, the primal iteration \eqref{eq:x:update} minimizes the augmented Lagrangian plus a proximal term $\frac{\beta}{2}\|x-x^r\|^2_{ B^T B} $. It is important to note that the proximal term is critical  in both the algorithm implementation and the analysis. It is used to ensure the following  key properties:  
\begin{enumerate}
	\item The primal problem is strongly convex, hence easily solvable;
	\item The primal problem is decomposable over different variable blocks.
\end{enumerate}
To see why the first point above is possible, suppose $B^T B$ is chosen such that $A^T A + B^T B\succeq I$, and that $f(x)$ has Lipschitz gradient. Then by a result in \cite[Theorem 2.1]{ZLOBEC05}, we know that for any $\beta>L$, the objective function of the $x$-problem \eqref{eq:x:update} is strongly convex. 

We illustrate the second point through an example. Consider the distributed optimization problem \eqref{eq:global:consensus:equiv}. Define the {\it signless incidence matrix} $B:=|A|$,  where the absolute value is taken for each component of $A$, and $A$ is the singed incidence matrix defined in Section \ref{sec:introduction}. Using this choice of $B$, we have $B^T B = L_{+}\in\mathbb{R}^{N\times N}$, which is the signless graph Laplacian whose $(i,i)$th diagonal entry is the degree of node $i$, and its $(i,j)$th entry is $1$ if $e=(i,j)\in\cE$, and $0$ otherwise.  Then $x$-update step \eqref{eq:x:update} becomes  
\begin{align*}
x^{r+1}& =\arg\min_{x}\; \sum_{i=1}^{N}f_i(x_i) +\langle\mu^r, Ax - b\rangle+\frac{\beta}{2}x^T L_{-} x+\frac{\beta}{2}(x-x^r)^T L_{+}(x-x^r)\\
& = \arg\min_{x}\;\sum_{i=1}^{N}f_i(x_i) +\langle\mu^r, Ax - b\rangle+\frac{\beta}{2}x^T (L_{-}+L_{+}) x-{\beta}x^T L_{+}x^{r}				\\
& = \arg\min_{x}\;\sum_{i=1}^{N}f_i(x_i) +\langle\mu^r, Ax - b\rangle+{\beta}x^T D x-{\beta}x^T L_{+}x^{r}				
\end{align*}
where $D=\mbox{diag}[d_1,\cdots, d_N]\in\mathbb{R}^{N\times N}$ is the degree matrix, with $d_i$ denoting the degree of node $i$. Clearly this problem is separable over the nodes, therefore it can be solved completely distributedly.  

We remark that one can always add an additional proximal term  $1/2\|x^{r}-x\|^2_{W}$ ($W\in\mathbb{R}^{N\times N}$ is some positive semidefinite matrix) to the $x$-subproblem \eqref{eq:x:update}. Our  analysis (to be presented shortly) remains valid. However in order to reduce the notational burden, we will solely focus on analyzing Algorithm 1 as presented in the table.  Also, Algorithm 1 can have many useful extensions. For example, one can solve the $x$-subproblem {\it inexactly}, by performing a proximal gradient step (in which $\langle \nabla f(x^r), x-x^r\rangle$ is used in place of $f(x)$ in \eqref{eq:x:update}). We will discuss these extensions in later sections.

\section{The Convergence Analysis}\label{sec:analysis}
In this section we provide convergence analysis for Algorithm 1. The key is the construction of a new potential function that decreases at every iteration of the algorithm. The constructed potential function is a conic combination the augmented Lagrangian, certain proximal term as well as the size of the violation of the equality constraint, thus it measures the progress of both the primal and dual updates.  

\subsection{The Assumptions}


We first state our main assumptions. 
\begin{itemize}
	\item[A1.] The function $f(x)$ is differentiable and has Lipschitz continuous gradient, i.e.,
	$$\|\nabla f(x) - \nabla f(y)\|\le L\|x-y\|, \quad \forall~x,y\in\mathbb{R}^N.$$
	Further assume that $A^T A + B^T B\succeq I_N.$
	\item[A2.] There exists a constant $\delta>0$ such that 
	$$\exists \underline{f}>-\infty, \quad \st\; f(x)+\frac{\delta}{2}\|Ax-b\|^2\ge \underline{f},\; \forall~x\in\mathbb{R}^N.$$
	
	{Without loss of generality and for the simplicity of notations, in the following we will assume that $\underline{f} =0$ \footnote{We note that this is without loss of generality because we have assumed that $\underline{f}$ is finite, therefore we can consider an {\it equivalent} problem with the objective function $f(x)+\underline{f}$, which is always lower bounded by $0$.}. }
	
	\item [A3.] The constraint $Ax = b$ is feasible over $x\in\mathbb{R}^N$.  
\end{itemize}
Below we provide a few nonconvex smooth functions $f(x)$ that satisfy Assumption [A1] -- [A3]. Note that the first three nonconvex functions are of particular interest in learning neural networks, as they are commonly used as  activation functions.
	\begin{itemize}
		\item {\bf The $\mbox{sigmoid}$ function}. The sigmoid function is given by
		$$\mbox{sigmoid}(x) = \frac{1}{1+e^{-x}}\in [-1,\;1].$$
		Clearly it satisfies [A2]. We have $\mbox{sigmoid}'(x) = \frac{e^{-x}}{(1+e^{-x})^2}\in [0, \; 1/4]$, and such boundedness of the first order derivative implies that [A1] is true (by first-order mean value theorem).  
		\item {\bf The $\arctan$ function}. Note that $\arctan(x)\in [-1, 1]$, so  it clearly satisfies [A2]. $\arctan'(x) = \frac{1}{x^2 +1} \in [0,\;1]$ so  it is bounded, which implies that [A1] is true. Finally note that 
		$$\arctan''(x) =  -2\frac{x}{(1+x^2)^2}$$
		so  it is a nonconvex function. 
		
		\item {\bf The $\tanh$ function}. Note that we have 
		$$\tanh(x)\in [-1, 1], \quad \tanh'(x) =1 - \tanh(x)^2 \in [0,1].$$ 
		Therefore the function satisfies [A1]--[A2]. 
		
		\item {\bf The} $\mbox{logit}$ {\bf function}. Since the logistic function is related to the $\tanh$ function as follows
		$$2\mbox{logit}(x) =\frac{2e^{x}}{e^x+1}=1+ \tanh(x/2),$$ then  Assumptions [A1]-[A2] are again satisfied.

		\item {\bf The } $\log(1+ x^2)$ {\bf function}.  This function has applications in structured matrix factorization \cite{fu16}. The function itself is obviously nonconvex and lower bounded. Its first order derivative is $\log'(1+ x^2) = \frac{2x}{1+x^2}$ and it is also bounded. 
		
		\item {\bf The quadratic function $x^T Q x$.} Suppose that $Q$ is a symmetric matrix but not necessarily positive semidefinite, and suppose that $x^T Q x$ is strongly convex in the null space of $A^T A$. Then it can be shown that there exists a $\delta$ large enough such that [A2] is true; see e.g., \cite{zhang10,bertsekas82}.
	\end{itemize}
Other relevant functions include $\sin(x)$, $\mbox{sinc}(x)$, $\cos(x)$ and so on. 

\subsection{The Analysis Steps}
Below we provide the analysis of Prox-PDA. Our analysis consists of a series of lemmas leading to a theorem characterizing the convergence and iteration complexity for prox-PDA. 

In the first step we provide a bound of the size of the constraint violation using a quantity related to the primal iterates. Let $\sigma_{\min}(A^T A)$ denote the smallest {\it non-zero} eigenvalue for $A^T A$. We have the following result.   
\begin{lemma}\label{lemma:mu:bound}
{\it 	Suppose Assumptions [A1] and [A3] are satisfied. Then the following is true for Prox-PDA. 
	\begin{align}\label{eq:mu:difference:bound}
	\frac{1}{\beta}\|\mu^{r+1}-\mu^{r}\|^2&\le \frac{2L^2}{\beta\sigma_{\min}(A^T A)}\left\| x^{r}-x^{r+1}\right\|^2\nonumber\\
	&\quad +\frac{2\beta}{\sigma_{\min}(A^T A)}\left\|B^T B \left((x^{r+1}-x^r) -(x^{r}-x^{r-1})\right)\right\|^2.
	\end{align}} 
\end{lemma} 

\noindent{\it Proof.} From the optimality condition of the $x$ problem \eqref{eq:x:update} we have
\begin{align*}
&\nabla f(x^{r+1})+ A^T \mu^r + \beta A^T (Ax^{r+1}-b)+\beta B^T B (x^{r+1}-x^r) =0.
\end{align*}
Applying \eqref{eq:mu:update}, we have
\begin{align}\label{eq:mu:expression}
&A^T \mu^{r+1} = - \nabla f(x^{r+1})-\beta B^T B (x^{r+1}-x^r).
\end{align}

Note that by Assumption [A3], we have that $b$ lies in the column space of $A$. Therefore we must have 
$$\mu^{r+1}-\mu^r = \beta (A x^{r+1}-b)\in \mbox{col}(A).$$
That is, the difference of the dual variable lies in the column space of $A$. 
Also from the fact that $\mu^0=0$, we have that the dual variable itself also lies in the column space of $A$
\begin{align}\label{eq:mu:column}
\mu^r = \beta\sum_{t=1}^{r} (A x^t-b)\in \mbox{col} (A), \quad \forall~r=1,\cdots.
\end{align}
Applying these facts to \eqref{eq:mu:expression}, 
and let $\sigma_{\min}(A^TA)$ denote the smallest {\it non-zero} eigenvalue of $A^T A$, we have
$$\sigma^{1/2}_{\min}(A^T A) \|\mu^{r+1}-\mu^r\|\le \|A(\mu^{r+1}-\mu^r)\|.$$
This inequality combined with \eqref{eq:mu:expression} implies that
\begin{align}
\|\mu^{r+1}-\mu^{r}\|&\le\frac{1}{\sigma^{1/2}_{\min}(A^T A)}\|- \nabla f(x^{r+1})-\beta B^T B (x^{r+1}-x^r) - (- \nabla f(x^{r})-\beta B^T B (x^{r}-x^{r-1}))\|\nonumber\\
& = \frac{1}{\sigma^{1/2}_{\min}(A^T A)}\left\|\nabla f(x^{r})- \nabla f(x^{r+1})-\beta B^T B \left((x^{r+1}-x^r) -(x^{r}-x^{r-1})\right)\right\|\nonumber.
\end{align}
Squaring both sides and dividing by $\beta$, we obtain the desired result. \QED

Our second step bounds the descent of the augmented Lagrangian. 
\begin{lemma}\label{lemma:lag:bound}
		{\it Suppose Assumptions [A1] and [A3] are satisfied. Then the following is true for Algorithm 1

\begin{align}\label{eq:lag:bound}
\begin{split}
L_{\beta}(x^{r+1},\mu^{r+1}) - L_{\beta}(x^r, \mu^r)
&\le -\left(\frac{\beta-L}{2}-\frac{2L^2}{\beta\sigma_{\min}(A^T A)}\right)\|x^{r+1}-x^r\|^2 \\
&\quad +\frac{2\beta\|B^TB\|}{\sigma_{\min}(A^T A)}\left\| \left((x^{r+1}-x^r) -(x^{r}-x^{r-1})\right)\right\|_{B^TB}^2.
\end{split}
\end{align}}
\end{lemma}
\noindent{\it Proof.} Since $f(x)$ has Lipschitz continuous gradient, and that $A^T A + B^T B\succeq I$ by Assumption [A1], it is known that if $\beta>L$, then the $x$-subproblem \eqref{eq:x:update} is strongly convex with modulus $\gamma: = \beta-L>0$; cf. \cite[Theorem 2.1]{ZLOBEC05}. That is, we have
\begin{align}\label{eq:strong}
&L_{\beta}(x, \mu^r) +\frac{\beta}{2} \|x-x^r\|^2_{B^T B}- (L_{\beta}(z, \mu^r)+\frac{\beta}{2} \|z-x^r\|^2_{B^T B})\nonumber\\
&\ge \langle \nabla_x L_{\beta}(z, \mu^r) + \beta (B^T B(z-x^r)), x-z\rangle +\frac{\gamma}{2}\|x-z\|^2, \;\forall~x,z\in\mathbb{R}^N, \; \forall~\mu^r.
\end{align}

Using this property, we have
\begin{align}\label{eq:lag:bound:derive}
&L_{\beta}(x^{r+1},\mu^{r+1}) - L_{\beta}(x^r, \mu^r)\nonumber\\
&=L_{\beta}(x^{r+1},\mu^{r+1}) - L_{\beta}(x^{r+1},\mu^{r}) + L_{\beta}(x^{r+1},\mu^{r})  - L_{\beta}(x^r, \mu^r)\nonumber\\
&\le L_{\beta}(x^{r+1},\mu^{r+1}) - L_{\beta}(x^{r+1},\mu^{r}) + L_{\beta}(x^{r+1},\mu^{r})  + \frac{\beta}{2} \|x^{r+1}-x^r\|^2_{B^T B}  - L_{\beta}(x^r, \mu^r)\nonumber\\
&\stackrel{\rm(i)}\le  \frac{\|\mu^{r+1}-\mu^r\|^2}{\beta}  +\langle\nabla_x L_{\beta}(x^{r+1}, \mu^r) + \beta (B^T B(x^{r+1}-x^r)), x^{r+1}-x^{r}\rangle -\frac{\gamma}{2}\|x^{r+1}-x^r\|^2 \nonumber\\
&\stackrel{\rm(ii)}\le -\frac{\gamma}{2}\|x^{r+1}-x^r\|^2 + \frac{\|\mu^{r+1}-\mu^r\|^2}{\beta}\nonumber\\
&\le -\frac{\gamma}{2}\|x^{r+1}-x^r\|^2 +  \frac{1}{\sigma_{\min}(A^T A)}\left(\frac{2L^2}{\beta}\left\| x^{r}-x^{r+1}\right\|^2+{2\beta}\left\|B^T B \left((x^{r+1}-x^r) -(x^{r}-x^{r-1})\right)\right\|^2\right)\nonumber\\
&= -\left(\frac{\beta-L}{2}-\frac{2L^2}{\beta\sigma_{\min}(A^T A)}\right)\|x^{r+1}-x^r\|^2 +\frac{2\beta}{\sigma_{\min}(A^T A)}\left\|B^T B \left((x^{r+1}-x^r) -(x^{r}-x^{r-1})\right)\right\|^2
\end{align}
where in $(i)$ we have used \eqref{eq:strong} with the identification $z=x^{r+1}$ and $x=x^r$; in $(ii)$ we have used the optimality condition for the $x$-subproblem \eqref{eq:x:update}.  The claim is proved. \QED

A key observation from Lemma \ref{lemma:lag:bound} is that no matter how large  $\beta$ is, the rhs of \eqref{eq:lag:bound} cannot be made negative, as the second term is {\it increasing} in $\beta$. This observation suggests that the augmented Lagrangian alone cannot serve as the potential function for Prox-PDA. 

In search for an appropriate potential function, we need a new object that is decreasing in the order of $\beta\left\|(x^{r+1}-x^r) -(x^{r}-x^{r-1})\right\|_{B^TB}^2$. The following lemma shows that the descent of the sum of the constraint violation $\|Ax^{r+1}-b\|^2$ and the proximal term $\|x^{r+1}-x^r\|_{B^T B}^2$ has the desired term.

\begin{lemma}\label{lemma:second:descent}
	Suppose Assumption [A1] is satisfied. Then the following is true
	\begin{align}\label{eq:key:bound}
	&\frac{\beta}{2}\left(\|Ax^{r+1}-b\|^2+\|x^{r+1}-x^r\|^2_{B^T B}\right)\nonumber\\
	&\le L\|x^{r+1}-x^r\|^2 +\frac{\beta}{2}\left(\|x^{r}-x^{r-1}\|^2_{B^T B}+ \|Ax^{r}-b\|^2\right)\nonumber\\
	&\quad -\frac{\beta}{2}\left(\|(x^{r}-x^{r-1})-(x^{r+1}-x^r)\|_{B^T B}^2+\|A(x^{r+1}-x^r)\|^2\right).
	\end{align}
\end{lemma}

\noindent{\it Proof.} From the optimality condition of the $x$-subproblem \eqref{eq:x:update} we have
\begin{align}
\langle \nabla f(x^{r+1})+A^T \mu^r+\beta A^T (A x^{r+1}-b) + \beta B^T B(x^{r+1}-x^r), x^{r+1}-x\rangle\le 0,\; \forall~x\in\mathbb{R}^N\nonumber\\
\langle \nabla f(x^{r})+A^T \mu^{r-1}+\beta A^T (A x^{r}-b) + \beta B^T B(x^{r}-x^{r-1}), x^{r}-x\rangle\le 0,\; \forall~x\in\mathbb{R}^N\nonumber.
\end{align}
Plugging $x=x^r$ into the first inequality and $x=x^{r+1}$ into the second, adding the resulting inequalities and utilizing the $\mu$-update step \eqref{eq:mu:update} we obtain
\begin{align}
\langle \nabla f(x^{r+1})-\nabla f(x^r)+A^T (\mu^{r+1}-\mu^r)+ \beta B^T B\left((x^{r+1}-x^r)-(x^r-x^{r-1})\right), x^{r+1}-x^r\rangle\le 0.\nonumber
\end{align}
Rearranging, we have
\begin{align}\label{eq:key}
&\langle  A^T (\mu^{r+1}-\mu^r), x^{r+1}-x^r\rangle\nonumber\\
&\le-\langle \nabla f(x^{r+1})-\nabla f(x^r)+ \beta B^T B\left((x^{r+1}-x^r)-(x^r-x^{r-1})\right), x^{r+1}-x^r\rangle .
\end{align} 
Let us bound the lhs and the rhs of \eqref{eq:key} separately. 

First the lhs of \eqref{eq:key} can be expressed as
\begin{align}\label{eq:key:lhs}
&\langle  A^T (\mu^{r+1}-\mu^r), x^{r+1}-x^r\rangle\nonumber\\
&= \langle  \beta A^T (Ax^{r+1}-b), x^{r+1}-x^r\rangle\nonumber\\
&= \langle  \beta (Ax^{r+1}-b), A x^{r+1}-b-(Ax^r-b)\rangle\nonumber\\
&= \beta \|Ax^{r+1}-b\|^2-\beta\langle A x^{r+1}-b, Ax^r-b\rangle\nonumber\\ 
& =\frac{\beta}{2}\left(\|Ax^{r+1}-b\|^2-\|Ax^{r}-b\|^2+\|A(x^{r+1}-x^r)\|^2\right).
\end{align}
Second we have the following bound for the rhs of \eqref{eq:key}
\begin{align}\label{eq:key:rhs}
&-\langle \nabla f(x^{r+1})-\nabla f(x^r)+ \beta B^T B\left((x^{r+1}-x^r)-(x^r-x^{r-1})\right), x^{r+1}-x^r\rangle \nonumber\\
&\le L\|x^{r+1}-x^r\|^2 -\beta \langle B^T B\left((x^{r+1}-x^r)-(x^r-x^{r-1})\right), x^{r+1}-x^r\rangle\nonumber\\ 
&=L\|x^{r+1}-x^r\|^2 +\frac{\beta}{2}\bigg(\|x^{r}-x^{r-1}\|^2_{B^T B}-\|x^{r+1}-x^r\|^2_{B^T B} -\|(x^{r}-x^{r-1})-(x^{r+1}-x^r)\|_{B^T B}^2\bigg).
\end{align}
Combining the above two bounds, we have
\begin{align*}
\frac{\beta}{2}\left(\|Ax^{r+1}-b\|^2+\|x^{r+1}-x^r\|^2_{B^T B}\right)&\le L\|x^{r+1}-x^r\|^2 +\frac{\beta}{2}\left(\|x^{r}-x^{r-1}\|^2_{B^T B}+ \|Ax^{r}-b\|^2\right)\nonumber\\
&\quad -\frac{\beta}{2}\left(\|(x^{r}-x^{r-1})-(x^{r+1}-x^r)\|_{B^T B}^2+\|A(x^{r+1}-x^r)\|^2\right).
\end{align*}
The desired claim is proved. \QED

It is interesting to observe that the new object, ${\beta}/{2}\left(\|Ax^{r+1}-b\|^2+\|x^{r+1}-x^r\|^2_{B^T B}\right)$, {\it increases} in $L\|x^{r+1}-x^r\|^2$ and {\it decreases} in $\frac{\beta}{2}\left(\|(x^{r}-x^{r-1})-(x^{r+1}-x^r)\|_{B^T B}^2\right)$, while the augmented Lagrangian behaves in an opposite manner (cf. Lemma \ref{lemma:lag:bound}). More importantly, in our new object, the constant in front of  $\|x^{r+1}-x^r\|^2$ is {\it independent} of $\beta$. Although neither of these two objects decreases by itself, quite surprisingly, a proper {\it conic combination} of these two objects decreases at every iteration of Prox-PDA. To precisely state the claim, let us define the {\it potential function} for Algorithm 1 as
\begin{align}\label{eq:potential}
P_{c, \beta}(x^{r+1}, x^r, \mu^{r+1}) = L_{\beta}(x^{r+1},\mu^{r+1}) + \frac{c\beta}{2}\left(\|Ax^{r+1}-b\|^2+\|x^{r+1}-x^r\|^2_{B^T B}\right)
\end{align}
where $c>0$ is some constant to be determined later. We have the following result. 
\begin{lemma}\label{lemma:potential:decrease}
{\it Suppose the assumptions made in Lemmas \ref{lemma:mu:bound} -- \ref{lemma:second:descent} are satisfied. 
Then we have the following  estimate
\begin{align}\label{eq:descent:potential}
P_{c,\beta}(x^{r+1}, x^r, \mu^{r+1})&\le P_{c,\beta}(x^{r}, x^{r-1}, \mu^{r})-\left(\frac{\beta-L}{2}-\frac{2L^2}{\beta\sigma_{\min}(A^T A)}-cL\right)\|x^{r+1}-x^r\|^2 \nonumber\\
&\quad -\left(\frac{c \beta}{2}-\frac{2\beta\|B^T B\|}{\sigma_{\min}(A^T A)}\right)\left\|(x^{r+1}-x^r) -(x^{r}-x^{r-1})\right\|_{B^T B}^2.
\end{align}}
\end{lemma}
\noindent {\it Proof}. Multiplying both sides of \eqref{eq:key:bound} by the constant $c$ and then add them to \eqref{eq:lag:bound}, we obtain
\begin{align*}
&L_{\beta}(x^{r+1},\mu^{r+1}) + \frac{c\beta}{2}\left(\|Ax^{r+1}-b\|^2+\|x^{r+1}-x^r\|^2_{B^T B}\right)\nonumber\\
&\le L_{\beta}(x^r, \mu^r)+ cL\|x^{r+1}-x^r\|^2 +\frac{c\beta}{2}\left(\|x^{r}-x^{r-1}\|^2_{B^T B}+ \|Ax^{r}-b\|^2\right)\nonumber\\
&-\left(\frac{\beta-L}{2}-\frac{2L^2}{\beta\sigma_{\min}(A^T A)}\right)\|x^{r+1}-x^r\|^2 +\frac{2\beta}{\sigma_{\min}(A^T A)}\left\|B^T B \left((x^{r+1}-x^r) -(x^{r}-x^{r-1})\right)\right\|^2\nonumber\\
&-\frac{c\beta}{2}\left(\|(x^{r}-x^{r-1})-(x^{r+1}-x^r)\|_{B^T B}^2+\|A(x^{r+1}-x^r)\|^2\right)\nonumber\\
&\le L_{\beta}(x^r, \mu^r)+ \frac{c\beta}{2}\left(\|x^{r}-x^{r-1}\|^2_{B^T B}+ \|Ax^{r}-b\|^2\right)\nonumber\\
&\quad -\left(\frac{\beta-L}{2}-\frac{2L^2}{\beta\sigma_{\min}(A^T A)}-cL\right)\|x^{r+1}-x^r\|^2 \nonumber\\
&\quad -\left(\frac{c \beta}{2}-\frac{2\beta\|B^T B\|}{\sigma_{\min}(A^T A)}\right)\left\|(x^{r+1}-x^r) -(x^{r}-x^{r-1})\right\|_{B^T B}^2.
\end{align*}
The desired result is proved.  \QED

From the above analysis, it is easy to see that as long as $c$ and $\beta$ are chosen large enough, the potential function decreases at each iteration of Prox-PDA. Below we derive the precise bounds for $c$ and $\beta$. 

First, it is clear that a sufficient condition for $c$ is 
\begin{align}
c\ge\max\left\{\frac{\delta}{L}, \frac{4\|B^T B\|}{\sigma_{\min}(A^T A)}\right\}\label{eq:c:bound}.
\end{align}
Note that the term``$\delta/L$" (Defined in Assumption [A2]) in the $\max$ operator is needed for later use.  Importantly, such bound on $c$ is {\it independent} of $\beta$.

Second, for any given $c$, we need $\beta$ to satisfy
\begin{align*}
\frac{\beta-L}{2}-\frac{2L^2}{\beta\sigma_{\min}(A^T A)}-cL>0,
\end{align*}
which implies the following
\begin{align}\label{eq:beta:bound}
\beta> \frac{L}{2}\left(2c+1+\sqrt{(2c+1)^2+\frac{16 L^2}{\sigma_{\min}(A^T A)}}\right).
\end{align}
Clearly combining the bounds for $\beta$ and $c$ we see that $\beta>\delta$. 
We conclude that if both \eqref{eq:c:bound} and \eqref{eq:beta:bound} are satisfied, then the potential function $P(x^{r+1}, x^r, \mu^{r+1})$ decreases at every iteration. 

Our next step shows that by using the particular choices of $c$ and $\beta$ in \eqref{eq:c:bound} and \eqref{eq:beta:bound}, the constructed potential function is lower bounded.  
\begin{lemma}\label{lemma:lower:bound}
{\it Suppose [A1] - [A3] are satisfied, and $(c,\beta)$ are chosen according to \eqref{eq:c:bound} and \eqref{eq:beta:bound}. Then the following statement holds true
\begin{align}\label{eq:lower:bound}
\exists~\underline{P}  \quad \st \;  P_{c,\beta}(x^{r+1}, x^r, \mu^{r+1}) \ge \underline{P}>-\infty, \quad \forall~r>0.
\end{align}}
\end{lemma}

\noindent{\it Proof}. To prove this we  need to utilize the boundedness assumption in [A2]. 

First, we can express the augmented Lagrangian function as following
\begin{align*}
L_{\beta}(x^{r+1}, \mu^{r+1})& = f(x^{r+1})+\langle \mu^{r+1}, A x^{r+1}-b\rangle +\frac{\beta}{2}\|A x^{r+1}-b\|^2\nonumber\\
&=f(x^{r+1})+\frac{1}{\beta}\langle \mu^{r+1}, \mu^{r+1}-\mu^r\rangle +\frac{\beta}{2}\|A x^{r+1}-b\|^2\nonumber\\
&=f(x^{r+1})+\frac{1}{2\beta}\left(\|\mu^{r+1}\|^2-\|\mu^r\|^2+\|\mu^{r+1}-\mu^r\|^2\right)+\frac{\beta}{2}\|A x^{r+1}-b\|^2.
\end{align*}
Therefore, summing over $r=1\cdots, T$, we obtain
\begin{align*}
\sum_{r=1}^{T}L_{\beta}(x^{r+1}, \mu^{r+1})&=\sum_{r=1}^{T}\left(f(x^{r+1})+\frac{\beta}{2}\|A x^{r+1}-b\|^2 + \frac{1}{2\beta}\|\mu^{r+1}-\mu^r\|^2\right) \nonumber\\
&\quad +\frac{1}{2\beta}\left(\|\mu^{T+1}\|^2 -\|\mu^1\|^2\right).
\end{align*}
Suppose Assumption [A2] is satisfied and $\beta$ is chosen according to \eqref{eq:beta:bound} and \eqref{eq:c:bound}, then clearly the above sum is lower bounded since  { $$f(x)+\frac{\beta}{2}\|A x-b\|^2\ge f(x)+\frac{\delta}{2}\|A x-b\|^2\ge 0, \; \forall~x\in\mathbb{R}^N.$$} 
This fact implies that the sum of the potential function is also lower bounded (note, the remaining terms in the potential function are all nonnegative), that is
\begin{align*}
\sum_{r=1}^{T}P_{c, \beta}(x^{r+1},x^r,\mu^{r+1}) >-\infty, \quad \forall~T>0.
\end{align*}
Note that if $c$ and $\beta$ are chosen according to \eqref{eq:c:bound} and \eqref{eq:beta:bound}, then $P_{c,\beta}(x^{r+1},x^r,\mu^{r+1}) $ is nonincreasing. Combined with the lower boundedness of the sum of the potential function, we can conclude that the following is true
\begin{align}\label{eq:lower:bound:case1}
P_{c,\beta}(x^{r+1},x^r,\mu^{r+1}) >-\infty, \quad \forall~r>0.
\end{align}

This completes the proof. \QED

Now we are ready to present the main result of this section on the convergence and the rate of convergence of Prox-PDA. To this end, 
define $Q(x^{r+1}, \mu^{r+1})$ as the optimality gap of problem \eqref{eq:original}, given by
\begin{align}\label{eq:Q}
&Q(x^{r+1}, \mu^{r}): = \|\nabla_{x}L_{\beta}(x^{r+1}, \mu^{r})\|^2+\|A x^{r+1}-b\|^2.
\end{align}
It is easy to see that $Q(x^{r+1}, \mu^{r}) \to 0$ implies that the limit point  $(x^*,\mu^*)$ is a KKT point of \eqref{eq:original} that satisfies the following conditions
\begin{align}\label{eq:stationary}
0=\nabla f(x^*)+ A^T \mu^*, \quad Ax^* =b.
\end{align}
To see this, we can first observe that $A x^{r+1}-b\to 0$, implying $\mu^{r+1}-\mu^r\to 0$, therefore the second condition in \eqref{eq:stationary} hold. Second,  using the fact that the first term in the optimality gap also goes to zero, we have
\begin{align*}
&\lim_{r\to\infty}\nabla f(x^{r+1}) + A^T \mu^{r}+\beta A^T (A x^{r+1}-b)  = \nabla f(x^*) + A^T \mu^* = 0.
\end{align*}
Therefore the first condition in \eqref{eq:stationary} is true.

In the following result we  show that the optimality gap $Q(\cdot)$ not only  decreases to zero, but does so in a sublinear manner. This is the main result of this section. 
\begin{theorem}\label{thm:main}
{\it Suppose Assumption A is satisfied. Further suppose that the conditions on $\beta$ and $c$ in \eqref{eq:beta:bound} and \eqref{eq:c:bound} are satisfied. Then we have the following claims for the sequence generated by Prox-PDA. \begin{itemize}
	\item {\bf (Eventual Feasibility).} The constraint is satisfied in the limit, i.e.,
	$$\lim_{r\to\infty} \mu^{r+1}-\mu^r\to 0, \; \lim_{r\to\infty }Ax^r \to b, \;\mbox{and }\;  \lim_{r\to\infty}x^{r+1}-x^r = 0.$$
	\item {\bf (Boundedness of Sequence)}. Further, if $\|\nabla f(x)\|$ is bounded for all $x\in\mathbb{R}^N$, then the sequence $\{\mu^r\}$ is also bounded;  If $f(x)+\frac{\beta}{2}\|Ax-b\|^2$ is coercive, then the sequence $\{x^r\}$ is bounded. 
	\item {\bf (Convergence to Stationary Points).} Every limit point of the iterates $\{x^r, \mu^r\}$ generated by Algorithm 1 converges to a stationary point of  problem \eqref{eq:original}.  Further, $Q(x^{r+1}, \mu^{r}) \to 0$. 
	\item {\bf (Sublinear Convergence Rate).} For any given $\varphi>0$, let us define $T$ to be the first time that the optimality gap reaches below $\varphi$, i.e., 
	$$T:=\arg\min_{r} Q(x^{r+1}, \mu^{r})\le \varphi.$$
	Then there exists a constant $\nu>0$ such that the following is true
	\begin{align*}
	\varphi\le \frac{\nu}{T-1}.
	\end{align*}
	That is, the optimality gap $Q(x^{r+1}, \mu^{r})$ converges sublinearly.
\end{itemize}}
\end{theorem}

\noindent {\it Proof.}  First we prove part (1). Combining Lemmas \ref{lemma:potential:decrease} and \ref{lemma:lower:bound}, we conclude that $\|x^{r+1}-x^r\|^2\to 0$. Then according to \eqref{eq:mu:difference:bound}, in the limit we have $\mu^{r+1}\to \mu^r$, or equivalently $Ax^r\to b$. That is, the constraint violation will be satisfied  in the limit.

Then we prove part (2). From the optimality condition of $x$-update step \eqref{eq:x:update} we have
\begin{align*}
&\nabla f(x^{r+1})+ A^T \mu^r + \beta A^T (Ax^{r+1}-b)+\beta B^T B (x^{r+1}-x^r) =0.
\end{align*}
Then we argue that $\{\mu^r\}$ is a bounded sequence if $\nabla f(x^{r+1})$ is bounded. Indeed  the fact that $\|x^{r+1}-x^r\|^2\to 0$ and $Ax^{r+1}\to b$ imply that both $(x^{r+1}-x^r)$ and $A x^{r+1}-b$ are bounded. Then the boundedness of $\mu^r$ follows from the assumption that $\nabla f(x)$ is bounded for any $x\in\mathbb{R}^N$, and that $\mu^r$ lies in the column space of $A$ (cf. \eqref{eq:mu:column}). 

Then we argue that $\{x^r\}$ is bounded if $f(x)+\frac{\beta}{2}\|Ax-b\|^2$ is coercive. Note that the potential function can be expressed as
\begin{align*}
P_{c,\beta}(x^{r+1}, x^r, \mu^{r+1}) &= f(x^{r+1})+\langle \mu^{r+1}, Ax^{r+1}-b\rangle+\frac{\beta}{2}\|Ax^{r+1}-b\|^2+ \frac{c\beta}{2}\left(\|Ax^{r+1}-b\|^2+\|x^{r+1}-x^r\|^2_{B^T B}\right)\nonumber\\
&= f(x^{r+1})+\frac{1}{2\beta}(\|\mu^{r+1}\|^2-\|\mu^r\|^2+\|\mu^{r+1}-\mu^r\|^2)+\frac{\beta}{2}\|Ax^{r+1}-b\|^2\nonumber\\
&\quad + \frac{c\beta}{2}\left(\|Ax^{r+1}-b\|^2+\|x^{r+1}-x^r\|^2_{B^T B}\right)
\end{align*}
and by our analysis in Lemma \ref{lemma:lower:bound} we know that it is decreasing thus {\it upper bounded}. 
Suppose that $\{x^r\}$ is unbounded and let $\cK$ denote an infinite subset of iteration index in which $\lim_{r\in \cK}x^r =\infty.$
Passing limit to $P_{c,\beta}(x^{r+1}, x^r, \mu^{r+1})$ over $\cK$, and using the fact that $x^{r+1}\to x^r$, $\mu^{r+1}\to \mu^r$, we have
\begin{align*}
\lim_{r\in\cK}P_{c,\beta}(x^{r+1}, x^r, \mu^{r+1}) = \lim_{r\in\cK} f(x^{r+1}) + \frac{c\beta+\beta}{2}\|Ax^{r+1}-b\| =\infty
\end{align*}
where the last equality comes from the coerciveness assumption. This is a contradiction to the fact that the potential function $P_{c,\beta}(x^{r+1}, x^r, \mu^{r+1}) $ is upper bounded. 
This concludes the proof for the second part of the result. 

Then we prove part (3). Let $\cal{K}$ denote any converging infinite iteration index such that $\{(\mu^r, x^r)\}_{r\in \cal{K}}$ converges to the limit point $(\mu^*, x^*)$. 
Passing limit in $\cal{K}$, and using the fact that $\|x^{r+1}-x^r\|\to 0$, we have
\begin{align*}
&\nabla f(x^*)+ A^T \mu^* + \beta A^T (Ax^{*}-b)=0.
\end{align*}
Combined with the fact that $Ax^*-b=0$, we conclude that $(\mu^*, x^*)$  is indeed a stationary point of the original problem \eqref{eq:original}, satisfying \eqref{eq:stationary}.

Additionally, even if the sequence $\{x^{r+1},\mu^{r+1}\}$ does not have a limit point, from part (1) we still have $\|\mu^{r+1}-\mu^r\|\to 0$ and $\|x^r-x^{r+1}\|\to 0$. Hence
\begin{align*}
\lim_{r\to\infty}\nabla_{x} L_{\beta}(x^{r+1}, \mu^{r})&=\lim_{r\to\infty}\nabla f(x^{r+1})+ A^T \mu^{r+1} +A^T (\mu^{r}-\mu^{r+1}) \nonumber\\
&\stackrel{\rm(i)}= \lim_{r\to\infty}  -\beta B^T B (x^{r+1}-x^r) + A^T (\mu^{r}-\mu^{r+1}) = 0
\end{align*}
where ${\rm(i)}$ is from the optimality condition of the $x$-subproblem \eqref{eq:x:update}. Therefore we have that $Q(x^{r+1}, \mu^{r})\to 0.$

Finally we prove part (4). Our first step is to bound the size of the gradient of the augmented Lagrangian. From the optimality condition of the $x$-problem \eqref{eq:x:update}, we have
\begin{align*}
\|\nabla_{x}L_{\beta}(x^{r}, \mu^{r-1})\|^2 &= \|\nabla_{x}L_{\beta}(x^{r+1}, \mu^{r})+\beta B^T B(x^{r+1}-x^r)- \nabla_{x}L_{\beta}(x^{r},\mu^{r-1})\|^2\nonumber\\
& = \|\nabla f(x^{r+1})-\nabla f(x^r) + A^T (\mu^{r+1}-\mu^{r}) +\beta B^T B (x^{r+1}-x^r)\|^2\nonumber\\
&\le 3 L^2\|x^{r+1}-x^r\|^2 + 3 \|\mu^{r+1}-\mu^r\|^2 \|A^T A\| + 3 \beta^2\|B^T B (x^{r+1}-x^r)\|^2.
\end{align*}
Therefore, by utilizing the estimate \eqref{eq:mu:difference:bound}, we see that there must exist two constants $\xi_1, \xi_2>0$ such that the following is true
\begin{align*}
Q(x^{r}, \mu^{r-1}) &= \|\nabla_{x}L_{\beta}(x^{r},\mu^{r-1})\|^2+\beta\|A x^{r}-b\|^2\nonumber\\
&\le\xi_1\left\| x^{r}-x^{r+1}\right\|^2+\xi_2\left\|B^T B \left((x^{r+1}-x^r) -(x^{r}-x^{r-1})\right)\right\|^2.
\end{align*}
From the descent estimate \eqref{eq:lag:bound} we see that there must exist two constants $\nu_1, \nu_2>0$ such that
\begin{align*}
&P_{c,\beta}(x^{r+1}, x^r, \mu^{r+1}) - P_{c,\beta}(x^{r}, x^{r-1}, \mu^{r})\nonumber\\
&\le -\nu_1\|x^{r+1}-x^r\|^2 - \nu_2\left\|B^T B \left((x^{r+1}-x^r) -(x^{r}-x^{r-1})\right)\right\|^2.
\end{align*}
Matching the above two bounds, we have
\begin{align*}
Q(x^{r}, \mu^{r-1}) \le \frac{\min\{\nu_1,\nu_2\}}{\max\{\xi_1,\xi_2\}} \left(P_{c,\beta}(x^{r}, x^{r-1}, \mu^{r}) - P_{c,\beta}(x^{r+1}, x^r, \mu^{r+1}) \right).
\end{align*}
Summing over $r$, and let $T$ denote the first time that $Q(x^{r+1}, x^r, \mu^{r+1})$ reaches below $\varphi$, we obtain
\begin{align*}
\varphi\le \frac{1}{T-1}\sum_{r=1}^{T}Q(x^{r}, \mu^{r-1}) &\le \frac{1}{T-1}\frac{\min\{\nu_1,\nu_2\}}{\max\{\xi_1,\xi_2\}}\left(P_{c,\beta}(x^{1}, x^{0}, \mu^{1}) - P_{c,\beta}(x^{T+1}, x^T, \mu^{T+1}) \right)\nonumber\\
&\le \frac{1}{T-1}\frac{\min\{\nu_1,\nu_2\}}{\max\{\xi_1,\xi_2\}}\left(P_{c,\beta}(x^{1}, x^{0}, \mu^{1}) - \underline{P}\right)\nonumber\\
&:= \frac{\nu}{T-1}.
\end{align*}
We conclude that the convergence in term of the optimality gap function  $Q(x^{r+1}, \mu^{r}) $  is sublinear.  \QED

Our result suggests that depending on the property of $f(x)$ and $\nabla f(x)$, the iterates $\{x^{r+1},\mu^{r+1}\}$ may or may not be bounded. However the optimality measure $Q(x,\mu)$ always converges to zero in a sublinear manner. Note that such sublinear complexity bound is in fact tight, even when applying first-order methods for nonconvex {\it unconstrained} problems; see the related discussions in \cite{Ghadimi14mini, Cartis10}.  

Before leaving this section, we remark that a few recent works \cite{hong14nonconvex_admm,hong14async_ADMM,li14nonconvex,Ames13LDA} have analyzed the convergence of a family of splitting method for certain nonconvex problems (which does not cover our formulation \eqref{eq:original}). All these works have used the augmented Lagrangian function as the potential function -- a technique first developed in \cite{Ames13LDA}. Unfortunately this technique fails to apply to our algorithm because it appears difficult, if not impossible, to show that the augmented Lagrangian alone achieves the desired descent (cf. Lemma \ref{lemma:lag:bound}).

\section{Discussion: The Convex Case}

It is interesting to observe that the proof techniques used in the previous section apply to the convex case as well -- only that for the convex case much milder conditions are required. That is, besides positivity, no additional requirement is needed  for the penalty parameter $\beta$. This observation also suggests that the proof techniques used here are rather ``tight", in the sense that it would be difficult to further sharpen the bounds on $\beta$ for the nonconvex case. 
It is interesting to observe that for convex cases Prox-PDA is closely related to the Method of Multipliers \cite{bertsekas82}, except that a proximal term is used in the $x$-step \eqref{eq:x:update}. Our analysis below shows that this type of method converges sublinearly for convex problems, without taking any averaging on the iterates, and for arbitrary choice of the positive penalty parameter. 

Below we briefly highlight the proof steps of the convex case.  Throughout this section, we will assume that Assumption [A1]-[A3] hold true. Additionally  assume that $f(x)$ is convex. 

First it is easy to bound the difference of the dual variables as 
\begin{align}
\|\mu^{r+1}-\mu^{r}\|&\le\frac{1}{\sigma^{1/2}_{\min}(A^T A)}\|- \nabla f(x^{r+1})-\beta B^T B (x^{r+1}-x^r) - (- \nabla f(x^{r})-\beta B^T B (x^{r}-x^{r-1}))\|\nonumber\\
& = \frac{1}{\sigma^{1/2}_{\min}(A^T A)}\left(\left\|\nabla f(x^{r})- \nabla f(x^{r+1})\right\|+\beta \left\|B^T B \left((x^{r+1}-x^r) -(x^{r}-x^{r-1})\right)\right\|\right)\nonumber.
\end{align}
The key difference compared with the proof in Lemma \ref{lemma:mu:bound} is that we do not use the Lipschitz continuity of $\nabla f(x)$ to make the rhs explicitly dependent on $\|x^{r+1}-x^r\|$.

Similarly as in Lemma \ref{lemma:lag:bound}, the descent of the augmented Lagrangian can be bounded by{
\begin{align*}
&L_{\beta}(x^{r+1},\mu^{r+1}) - L_{\beta}(x^r, \mu^r)\nonumber\\
&\le -\frac{\beta}{2}\|x^{r+1}-x^r\|^2 + \frac{\|\mu^{r+1}-\mu^r\|^2}{\beta}\nonumber\\
&\le -\frac{\beta}{2}\|x^{r+1}-x^r\|^2 +  \frac{2}{\beta\sigma_{\min}(A^T A)}\left\| \nabla f(x^{r})-\nabla f(x^{r+1})\right\|^2\nonumber\\
&\quad +\frac{2\beta}{\sigma_{\min}(A^T A)}\left\|B^T B \left((x^{r+1}-x^r) -(x^{r}-x^{r-1})\right)\right\|^2\nonumber.
\end{align*}}
\!\!Note that in the first inequality we have replaced $\gamma$ with $\beta$ due to the convexity assumption of $f(x)$ as well as the assumption that $A^T A + B^T B \succeq I$. 

Second, from Lemma \ref{lemma:second:descent} we obtain 
\begin{align}\label{eq:descent:second:convex}
&\frac{\beta}{2}\left(\|Ax^{r+1}-b\|^2+\|x^{r+1}-x^r\|^2_{B^T B}\right)\nonumber\\
&\le - \langle \nabla f(x^{r+1})-\nabla f(x^r), x^{r+1}-x^{r}\rangle  +\frac{\beta}{2}\left(\|x^{r}-x^{r-1}\|^2_{B^T B}+ \|Ax^{r}-b\|^2\right)\nonumber\\
&\quad -\frac{\beta}{2}\left(\|(x^{r}-x^{r-1})-(x^{r+1}-x^r)\|_{B^T B}^2+\|A(x^{r+1}-x^r)\|^2\right)\nonumber\\
&\le - \frac{1}{L}\|\nabla f(x^{r+1})-\nabla f(x^r)\|^2 +\frac{\beta}{2}\left(\|x^{r}-x^{r-1}\|^2_{B^T B}+ \|Ax^{r}-b\|^2\right)\nonumber\\
&\quad -\frac{\beta}{2}\left(\|(x^{r}-x^{r-1})-(x^{r+1}-x^r)\|_{B^T B}^2+\|A(x^{r+1}-x^r)\|^2\right)
\end{align}
where the last inequality is true due to the fact that the Lipschitz continuity of $\nabla f(x)$ and the convexity of $f(x)$ implies \cite{nesterov05}
\begin{align}\label{eq:nesterov}
\langle \nabla f(x)-\nabla f(z), x-z\rangle \ge \frac{1}{L}\|\nabla f(x)-\nabla f(z)\|^2,\; \forall~x,z\in\mathbb{R}^N.
\end{align}

Again define the potential function as 
\begin{align*}
P_{c, \beta}(x^{r+1}, x^r, \mu^{r+1}) = L_{\beta}(x^{r+1},\mu^{r+1}) + \frac{c\beta}{2}\left(\|Ax^{r+1}-b\|^2+\|x^{r+1}-x^r\|^2_{B^T B}\right).
\end{align*}
For any fixed $\beta>0$ pick  $c$ such that 
\begin{align}\label{eq:c:bound:convex}
c>\max\left\{\frac{2 L}{\beta \sigma_{\min}(A^T A)}, \frac{4 \|B^T B\|}{\sigma_{\min}(A^T A)}, \frac{\delta}{\beta}\right\}
\end{align}
then the potential function will decrease at every iteration of the algorithm, where the descent quantity is composed of a conic combination of the following terms: $-\|x^{r+1}-x^r\|^2$, $-\|\nabla f(x^r) - \nabla f(x^{r+1})\|$,  $-\|(x^{r}-x^{r-1})-(x^{r+1}-x^r)\|_{B^T B}^2$ and $-\|A (x^{r+1} - x^r)\|^2$.

The rest of the proof of the (rate of) convergence follows the similar arguments as those leading to Theorem \ref{thm:main}. To summarize this section, we provide the following corollary. 
\begin{corollary}\label{cor:convex}
	{\it Suppose Assumption [A1] - [A3] are satisfied. Suppose that $f(x)$ is convex, and $\beta$ is any positive number. Then the same conclusions in Theorem  \ref{thm:main} hold true for the Prox-PDA. }
\end{corollary}

\section{Extension: Inexactly Solving the Primal Problems}\label{sec:inexact}

In this section, we discuss two important extensions of the Prox-PDA, in which the $x$-problem \eqref{eq:x:update} is solved inexactly. Our first extension replaces the $x$-step by a single gradient-type step, while our second algorithm solves the $x$ problem \eqref{eq:x:update} to some predefined error. The motivation is that for many practical applications, exactly minimizing the augmented Lagrangian may not be easy. 

The proposed proximal gradient primal dual algorithm (Prox-GPDA) replaces the objective function $f(x)$ by the the surrogate function 
\begin{align}\label{eq:u}
u(x_1,x_2) :=\langle \nabla g(x_2), x_1-x_2\rangle .
\end{align}
The detailed algorithm is given in the following table. 
\begin{center}
	\fbox{
		\begin{minipage}{6in}
			\smallskip
			\centerline{\bf {Algorithm 2. The Proximal Gradient Primal Dual Algorithm (Prox-GPDA)}}
			\smallskip
			Initialize $\mu^0=0$ and  $x^0\in\mathbb{R}^N$;
			
			At each iteration $r+1$, update variables  by:
			\begin{subequations}
				\begin{align}
				x^{r+1}& =\arg\min_{x\in\mathbb{R}^N}\; u(x, x^r) +\langle\mu^r, Ax - b\rangle+\frac{\beta}{2}\|Ax-b\|^2+\frac{\beta}{2}(x-x^r)^T B^T B (x-x^r);\label{eq:x:update:1}\\
				\mu^{r+1}& = \mu^r +\beta (A x^{r+1}	-b)	.		\label{eq:mu:update:1}
				\end{align}
			\end{subequations}
			
		\end{minipage}
	}
\end{center}
\vspace{0.2cm}

Our second extension solves the $x$-step \eqref{eq:x:update} to certain $\epsilon$-optimality, where $\epsilon$ is some error term with small magnitude (the precise condition will be presented shortly). The quality of the solution is measured by the size of the gradient of the objective function for problem \eqref{eq:x:update}.

\begin{center}
	\fbox{
		\begin{minipage}{5.8in}
			\smallskip
			\centerline{\bf {Algorithm 3. The Inexact Proximal Primal Dual Algorithm (In-Prox-PDA)}}
			\smallskip
			Initialize $\mu^0$ and  $x^0$;
			
			At each iteration $r+1$, update variables  by:
			\begin{subequations}
				\begin{align}
				&\mbox{Find}\; x^{r+1} \; \mbox{that satisfies the following}\nonumber\\
				&\nabla f(x^{r+1})+ A^T \mu^r + \beta A^T(A x^{r+1}-b)+\beta B^T B(x^{r+1}-x^r) =\epsilon^{r+1}; \label{eq:x:update:2}\\
				&\mu^{r+1} = \mu^r +\beta (A x^{r+1}	-b)			\label{eq:mu:update:2}.
				\end{align}
			\end{subequations}
			
		\end{minipage}
	}
\end{center}
\vspace{0.2cm}

The analysis of these two cases follows similar steps as that for Prox-PDA. For Prox-GPDA, the major difference is that there are several places in which we need to bound the term $\|\nabla f(x^{r-1}) -\nabla f(x^r)\|$ instead of  $\|\nabla f(x^{r+1}) -\nabla f(x^r)\|$. Moreover, the potential function is no longer decreasing at each iteration. For the In-Prox-PDA case, an explicit condition on the size of the error sequence $\{\epsilon^{r+1}\}$ is needed. In the next subsection we provide an outline of the proof.

\subsection{The Analysis Outline for Prox-GPDA}

First, following the derivation leading to \eqref{eq:mu:difference:bound} we obtain
\begin{align}\label{eq:mu:difference:bound:1}
\frac{1}{\beta}\|\mu^{r+1}-\mu^{r}\|^2\le \frac{2L^2}{\beta\sigma_{\min}(A^T A)}\left\| x^{r}-x^{r-1}\right\|^2+\frac{2\beta}{\sigma_{\min}(A^T A)}\left\|B^T B \left((x^{r+1}-x^r) -(x^{r}-x^{r-1})\right)\right\|^2.
\end{align} 
Note that the first term is now related to the difference squared of the {\it previous} two iterations. 

Following the proof steps in Lemma \ref{lemma:lag:bound}, the descent of the augmented Lagrangian is given by
\begin{align}\label{eq:lag:bound:2}
&L_{\beta}(x^{r+1},\mu^{r+1}) - L_{\beta}(x^r, \mu^r)\nonumber\\
&\le -\frac{\beta-L}{2}\|x^{r+1}-x^r\|^2 +\frac{2\beta}{\sigma_{\min}(A^T A)}\left\|B^T B \left((x^{r+1}-x^r) -(x^{r}-x^{r-1})\right)\right\|^2\nonumber\\
&\quad + \frac{2L^2}{\beta\sigma_{\min}(A^T A)}\left\| x^{r}-x^{r-1}\right\|^2.
\end{align}

In the third step we have the following estimate
\begin{align}\label{eq:key:bound:1}
&\frac{\beta}{2}\left(\|Ax^{r+1}-b\|^2+\|x^{r+1}-x^r\|^2_{B^T B}\right)\nonumber\\
&\le \frac{L}{2}\|x^{r-1}-x^r\|^2 +\frac{L}{2}\|x^{r+1}-x^r\|^2 +\frac{\beta}{2}\left(\|x^{r}-x^{r-1}\|^2_{B^T B}+ \|Ax^{r}-b\|^2\right)\nonumber\\
&\quad -\frac{\beta}{2}\left(\|(x^{r}-x^{r-1})-(x^{r+1}-x^r)\|_{B^T B}^2+\|A(x^{r+1}-x^r)\|^2\right).
\end{align}

Note that the first two terms come from the following estimate
\begin{align*}
-\langle x^{r+1}-x^r, \nabla f(x^r)-\nabla f(x^{r-1})\rangle &\le \frac{L}{2}\|x^{r+1}-x^r\|^2+\frac{1}{2L}\|\nabla g(x^r)-\nabla g(x^{r-1})\|^2\nonumber\\
&\le  \frac{L}{2}\|x^{r+1}-x^r\|^2+\frac{L}{2}\|x^r-x^{r-1}\|^2.
\end{align*}

In the fourth step we have the following overall descent estimate
\begin{align}\label{eq:lag:bound:1}
&L_{\beta}(x^{r+1},\mu^{r+1}) + \frac{c\beta}{2}\left(\|Ax^{r+1}-b\|^2+\|x^{r+1}-x^r\|^2_{B^T B}\right)\nonumber\\
&\le L_{\beta}(x^r, \mu^r)+ \frac{c\beta}{2}\left(\|x^{r}-x^{r-1}\|^2_{B^T B}+ \|Ax^{r}-b\|^2\right)\nonumber\\
&\quad -\left(\frac{\beta-L}{2}+\frac{cL}{2}\right)\|x^{r+1}-x^r\|^2+\left(\frac{2L^2}{\beta\sigma_{\min}(A^T A)}+\frac{cL}{2}\right)\|x^{r-1}-x^r\|^2 \nonumber\\
&\quad -\left(\frac{c \beta}{2}-\frac{2\beta\|B^T B\|}{\sigma_{\min}(A^T A)}\right)\left\|(x^{r+1}-x^r) -(x^{r}-x^{r-1})\right\|_{B^T B}^2.
\end{align}
Note that there is a slight difference between this descent estimate and our previous estimate \eqref{eq:descent:potential}, because now there is a positive term in the rhs, which involves $\|x^{r}-x^{r-1}\|^2$. Therefore the potential function is difficult to decrease by itself. Fortunately, such extra term can be bounded by the descent of the {\it previous} iteration. 
We can take the summation over all the iterations and obtain
\begin{align}
&L_{\beta}(x^{T+1},\mu^{T+1}) + \frac{c\beta}{2}\left(\|Ax^{T+1}-b\|^2+\|x^{T+1}-x^{T}\|^2_{B^T B}\right)\nonumber\\
&\le L_{\beta}(x^1, \mu^1)+ \frac{c\beta}{2}\left(\|x^{1}-x^{0}\|^2_{B^T B}+ \|Ax^{1}-b\|^2\right)+\left(\frac{2L^2}{\beta\sigma_{\min}(A^T A)}+cL\right)\|x^{0}-x^1\|^2\nonumber\\
&\quad -\sum_{r=1}^{T-1}\left(\frac{\beta-L}{2}-\frac{2L^2}{\beta\sigma_{\min}(A^T A)}-cL\right)\|x^{r+1}-x^r\|^2 \nonumber\\
&\quad -\sum_{r=1}^{T}\left(\frac{c \beta}{2}-\frac{2\beta\|B^T B\|}{\sigma_{\min}(A^T A)}\right)\left\|(x^{r+1}-x^r) -(x^{r}-x^{r-1})\right\|_{B^T B}^2\nonumber.
\end{align}
Clearly as long as the potential function is lower bounded, we have $x^{r+1}\to x^r$ and $x^{r+1}-x^r\to x^{r}-x^{r-1}$.
 The rest of the proof follows similar steps leading to Theorem \ref{thm:main}, hence is omitted. 

\begin{corollary}\label{cor:proximal}
	{\it Suppose Assumption [A1]--[A3] are satisfied. Suppose $\beta$ and $c$ satisfy \eqref{eq:beta:bound} and \eqref{eq:c:bound}. Then all the conclusions in Theorem  \ref{thm:main} hold true for the Prox-GPDA. }
\end{corollary}

We remark that we can replace the approximation function $u(x, x^r)$ by a larger family of {\it upper-bound} functions, following  the BSUM (Block Successive Upper Bound Minimization) framework \cite{Razaviyayn12SUM,hong15busmm_spm}. The analysis follows similar lines of argument presented in this section. 

\subsection{The Analysis Outline for In-Prox-PDA}
Similarly as in Lemma \ref{lemma:mu:bound}, we have the following bound
	\begin{align*}
	\frac{1}{\beta}\|\mu^{r+1}-\mu^{r}\|^2&\le \frac{3L^2}{\beta\sigma_{\min}(A^T A)}\left\| x^{r}-x^{r+1}\right\|^2+ \frac{3}{\beta}\|\epsilon^{r+1}-\epsilon^r\|^2\nonumber\\
	&\quad +\frac{3\beta}{\sigma_{\min}(A^T A)}\left\|B^T B \left((x^{r+1}-x^r) -(x^{r}-x^{r-1})\right)\right\|^2.
	\end{align*} 

Then following the proof steps in Lemma \ref{lemma:lag:bound}, the descent of the augmented Lagrangian is given by
\begin{align*}
&L_{\beta}(x^{r+1},\mu^{r+1}) - L_{\beta}(x^r, \mu^r)\nonumber\\
&\le -\frac{\beta-L}{2}\|x^{r+1}-x^r\|^2 +\frac{3\beta}{\sigma_{\min}(A^T A)}\left\|B^T B \left((x^{r+1}-x^r) -(x^{r}-x^{r-1})\right)\right\|^2\nonumber\\
&\quad + \frac{3L^2}{\beta\sigma_{\min}(A^T A)}\left\| x^{r}-x^{r+1}\right\|^2 + \frac{3}{\beta}\|\epsilon^{r+1}-\epsilon^r\|^2+ \langle \epsilon^{r+1}, x^{r+1}-x^r\rangle.
\end{align*}
The inexactness update results in two additional terms in the descent of the augmented Lagrangian.

In the third step we have the following estimate
\begin{align*}
&\frac{\beta}{2}\left(\|Ax^{r+1}-b\|^2+\|x^{r+1}-x^r\|^2_{B^T B}\right)\nonumber\\
&\le L\|x^{r+1}-x^r\|^2 +\frac{\beta}{2}\left(\|x^{r}-x^{r-1}\|^2_{B^T B}+ \|Ax^{r}-b\|^2\right)\nonumber\\
&\quad -\frac{\beta}{2}\left(\|(x^{r}-x^{r-1})-(x^{r+1}-x^r)\|_{B^T B}^2+\|A(x^{r+1}-x^r)\|^2\right) + \langle \epsilon^{r+1}-\epsilon^r, x^{r+1}-x^r\rangle .
\end{align*}
In the fourth step we have the following overall descent estimate
\begin{align*}
&L_{\beta}(x^{r+1},\mu^{r+1}) + \frac{c\beta}{2}\left(\|Ax^{r+1}-b\|^2+\|x^{r+1}-x^r\|^2_{B^T B}\right)\nonumber\\
&\le L_{\beta}(x^r, \mu^r)+ \frac{c\beta}{2}\left(\|x^{r}-x^{r-1}\|^2_{B^T B}+ \|Ax^{r}-b\|^2\right)\nonumber\\
&\quad -\left(\frac{\beta-L}{2}-\frac{3L^2}{\beta\sigma_{\min}(A^T A)}-cL-\frac{1+2c}{2}\right)\|x^{r+1}-x^r\|^2 \nonumber\\
&\quad -\left(\frac{c \beta}{2}-\frac{3\beta\|B^T B\|}{\sigma_{\min}(A^T A)}\right)\left\|\left((x^{r+1}-x^r) -(x^{r}-x^{r-1})\right)\right\|_{B^T B}^2\nonumber\\
&\quad +\frac{1}{2}\|\epsilon^{r+1}\|^2+\left(\frac{c}{4}+\frac{3}{\beta}\right)\|\epsilon^{r+1}-\epsilon^r\|^2.
\end{align*}
Therefore as long as the potential function is lower bounded, and that the error sequence satisfies 
\begin{align}\label{eq:error}
\sum_{r=0}^{\infty}\|\epsilon^{r+1}-\epsilon^r\|^2 <\infty, \quad \sum_{r=0}^{\infty}\|\epsilon^{r}\|^2 <\infty
\end{align}
we have $x^{r+1}\to x^r$ and $x^{r+1}-x^r\to x^{r}-x^{r-1}$.
The rest of the proof follows similar steps leading to Theorem \ref{thm:main}, hence are omitted. We have the following convergence for In-Prox-GPDA.

\begin{corollary}
	{\it Suppose Assumption [A1]--[A3] are satisfied. Suppose $\{\epsilon^r\}$ satisfies \eqref{eq:error}, and  $\beta$ and $c$ satisfy the following conditions
		\begin{align*}
		c&\ge\max\left\{\frac{\delta}{L}, \frac{6\|B^T B\|}{\sigma_{\min}(A^T A)}\right\},\nonumber\\
		\beta&> \frac{L+1}{2}\left((2c+1)+\sqrt{(2c+1)^2+\frac{16( L+1)^2}{\sigma_{\min}(A^T A)}}\right).
		\end{align*}
		Then all conclusions in Theorem  \ref{thm:main} hold true for the In-Prox-GPDA. }
\end{corollary}

\section{Extension: Increasing the Penalty Sequence}

In this section, we present an important variant of Prox-PDA in which there is no need to explicitly compute the bound for the penalty parameter $\beta$. Indeed, the bounds on $\beta$ derived in the previous sections are the worst case bounds, and algorithms that use stepsizes that strictly satisfy such bounds may be slow at the beginning.  In practice, one may prefer to start with a small penalty parameter and gradually increase it. The following algorithm adopts such strategy. 
\begin{center}
	\fbox{
		\begin{minipage}{5.6in}
			\smallskip
			\centerline{\bf {Algorithm 4. The Prox-PDA with Increasing Penalty (Prox-PDA-IP)}}
			\smallskip
			At iteration $0$, initialize $\mu^0 =0$ and $x^0\in\mathbb{R}^N$.
			
			At each iteration $r+1$, update variables  by:
			\begin{subequations}
				\begin{align}
				x^{r+1}& =\arg\min_{x\in\mathbb{R}^n}\; f(x) +\langle\mu^r, Ax - b\rangle+\frac{\beta^{r+1}}{2}\|Ax-b\|^2+\frac{\beta^{r+1}}{2}\|x-x^r\|^2_{B^T B}; \label{eq:x:update:inc}\\
				\mu^{r+1}& = \mu^r +\beta^{r+1} (A x^{r+1}	- b).		\label{eq:mu:update:inc}
				\end{align}
			\end{subequations}
			
		\end{minipage}
	}
\end{center}

We note that in the above algorithm, the primal proximal parameter, the primal penalty parameter, as well as the dual stepsize are all time-varying (in fact all of them increase unboundedly). This is the key feature of this variant. It would be challenging to achieve convergence if only a subset of these parameters grow unboundedly. 


Throughout this section we will still assume that Assumption A holds true. Further, we will assume that ${\beta^r}$ satisfies the following conditions
\begin{align}\label{eq:beta}
\frac{1}{\beta^r}\to 0, \quad \sum_{r=1}^{\infty}\frac{1}{\beta^r}=\infty, \quad \beta^{r+1}\ge \beta^r, \quad \max_{r}{(\beta^{r+1}-\beta^r)}<\omega, ~\mbox{for some finite  }~\omega>0.
\end{align} 
Also without loss of generality we will assume that 
\begin{align}\label{eq:B:inc}
B^T B\succ 0, \quad \mbox{and}\quad \|B^T B\|>1. 
\end{align}
Note that this is always possible, by adding an identity matrix to $B^T B$ if necessary.


The proof of convergence is long and technical, therefore we delegate it to Section \ref{sec:proof:increasing}. The key step is to construct a new potential function, given below
\begin{align}
P_{\beta^{r+1}, c}(x^{r+1}, x^r, \mu^{r+1})= L_{\beta^{r+1}}(x^{r+1},\mu^{r+1}) + \frac{c\beta^{r+1}\beta^r}{2}\|Ax^{r+1}-b\|^2 + \frac{c\beta^{r+1}\beta^r}{2}\|x^{r}-x^{r+1}\|^2_{B^T B}.\nonumber
\end{align}
Note that the key to construct the potential function for Algorithm 4 is to make the coefficients in front of the terms $\|x^{r+1}-x^r\|^2_{B^T B}$ and $\|Ax^{r+1}-b\|^2$ proportional to $(\beta^r)^2$. Our proof shows that after some finite number of iterations, the newly constructed potential function starts to descend, and the size of the descent is proportional to the following 
\begin{align}
\frac{\beta^{r+1}}{2}\|x^{r+1}-x^r\|^2+\frac{(\beta^{r})^2}{2}\|(x^{r+1}-x^r)-(x^r-x^{r-1})\|^2.
\end{align}

Therefore if the potential function is lower bounded, then we can conclude that 
\begin{align}
&\sum_{r=1}^{\infty}\; \frac{\beta^{r+1}}{2}\|x^{r+1}-x^r\|^2\le \infty\nonumber\\
&\sum_{r=1}^{\infty}\; \frac{(\beta^{r})^2}{2}\|(x^{r+1}-x^r)-(x^r-x^{r-1})\|^2\le \infty\nonumber. 
\end{align}
Using these two inequalities, we can show the desired convergence to stationary solution of problem \eqref{eq:original}. We refer the readers to to Section \ref{sec:proof:increasing} for proof details. 

We have the following theorem regarding to the convergence of Prox-PDA-IP. 
\begin{theorem}\label{thm:inc}
	{\it Suppose Assumption A is satisfied. Further suppose that the sequence of penalty parameters $\{\beta^r\}$ satisfies \eqref{eq:beta}, and that $B$ is selected such that \eqref{eq:B:inc} holds true. Then we have the following claims for Prox-PDA-IP. \begin{itemize}
			\item {\bf (Eventual Feasibility).} The constraint is satisfied in the limit, i.e.,
			$$\lim_{r\to\infty} \mu^{r+1}-\mu^r\to 0, \; \lim_{r\to\infty }Ax^r \to b, \;\mbox{and }\;  \lim_{r\to\infty}x^{r+1}-x^r = 0.$$
			\item {\bf (Boundedness of Sequence)}. Further, if $\|\nabla f(x)\|$ is bounded for all $x\in\mathbb{R}^N$, then the sequence $\{\mu^r\}$ is also bounded;  If there exists a constant $\theta>0$ such that $f(x)+\frac{\theta}{2}\|Ax-b\|^2$ is coercive, then the sequence $\{x^r\}$ is bounded. 
			\item {\bf (Convergence to Stationary Points).} Every limit point of the iterates $\{x^r, \mu^r\}$ generated by Algorithm 4 converges to a stationary point of  problem \eqref{eq:original}.  Further, $Q(x^{r+1}, \mu^{r}) \to 0$. 
		\end{itemize}}
	\end{theorem}

We remark the same analysis technique also applies to the following schemes, where $\beta^{r+1}$ also satisfies the condition in \eqref{eq:beta:bound}, and  the function $u(x,z)$ is given in \eqref{eq:u}.

\begin{center}
	\fbox{
		\begin{minipage}{6in}
			\smallskip
			\centerline{\bf {Algorithm 5. The Prox-GPDA with Increasing Penalty (Prox-GPDA-IP)}}
			\smallskip
			At iteration $0$, initialize $\mu^0 =0$ and $x^0\in\mathbb{R}^N$.
			
			At each iteration $r+1$, update variables  by:
			\begin{subequations}
				\begin{align}
				x^{r+1}& =\arg\min_{x\in\mathbb{R}^n}\; u(x, x^r) +\langle\mu^r, Ax-b\rangle+\frac{\beta^{r+1}}{2}\|Ax-b\|^2+\frac{\beta^{r+1}}{2}\|x-x^r\|^2_{B^T B}; \label{eq:x:update:inc:2}\\
				\mu^{r+1}& = \mu^r +\beta^{r+1} (A x^{r+1}	- b).			\label{eq:mu:update:inc:2}
				\end{align}
			\end{subequations}
			
		\end{minipage}
	}
\end{center}

%
%

\section{Application: Distributed Nonconvex Optimization}\label{sec:distributed}

In this section, we discuss the applications of Algorithms 1, 2 and 4 to the nonconvex distributed optimization problem. Our focus will be given to providing explicit update rules for each distributed node, as well as to relating the resulting algorithms to some well-known algorithms in the literature for distributed {\it convex} optimization.

We  assume throughout this section that each component function $f_i$ has Lipschitz continuous gradient with constant $L_i$. Then clearly Assumption [A1] is still satisfied, with $L=\max_i L_i$. 

First, we present a direct application of Prox-GPDA (Algorithm 2) to the distributed optimization setting. Quite interestingly, in this case we recover the so-called  EXTRA algorithm \cite{shi14extra} for distributed {\it convex} optimization. 

The derivation is in fact rather straightforward. The optimality condition of the $x$-update step \eqref{eq:x:update:1} is given by
\begin{align*}
\nabla f(x^r) +A^T \mu^r +\beta A^T Ax^{r+1} +\beta (B^T B(x^{r+1}-x^r))=0.
\end{align*} 
Utilizing the fact that $A^T A = L_{-}$, $B^T B = L_{+}$ and $L_{+}+L_{-}= 2D$ (where $L_{-}$, $L_{+}$ and $D$ denote respectively the signed Laplacian, the signless Laplacian and the degree matrix), we have
\begin{align*}
\nabla f(x^r) +A^T \mu^r +\beta 2 D x^{r+1} - \beta L_{+}x^r=0.
\end{align*} 
Subtracting the same equation evaluated at the previous iteration, we obtain
\begin{align*}
\nabla f(x^r) - \nabla f(x^{r-1}) + \beta L_{-} x^r +\beta 2 D (x^{r+1}-x^r) - \beta L_{+}(x^{r}-x^{r-1}) = 0
\end{align*}
where we have used the fact that $A^T (\mu^{r}-\mu^{r-1})=\beta A^T A x^r = \beta L_{-}x^r$. Rearranging terms, we have
\begin{align}\label{eq:EXTRA}
x^{r+1} &= x^r -\frac{1}{2\beta} D^{-1}\left(\nabla f(x^r) - \nabla f(x^{r-1})\right) + \frac{1}{2} D^{-1} (L_{+}-L_{-})x^r -\frac{1}{2} D^{-1} L_{+}x^{r-1}\nonumber\\
& = x^r -\frac{1}{2\beta} D^{-1}\left(\nabla f(x^r) - \nabla f(x^{r-1})\right) + W x^r -\frac{1}{2}(I+W)x^{r-1}
\end{align}
where in the last equality we have defined the {\it weight matrix} $W: = \frac{1}{2}D^{-1} (L_{+}-L_{-})$, which is a row stochastic matrix. 

Iteration \eqref{eq:EXTRA} has exactly the same form as the EXTRA algorithm given in \cite{shi14extra}, therefore we can conclude that the EXTRA is a special case of Prox-GPDA.  Moreover, by appealing to our analysis in Section \ref{sec:inexact}, it readily follows that iteration \eqref{eq:EXTRA} works for the nonconvex distributed optimization problem as well, as long as the parameter $\beta$ is selected appropriately. In particular, in this setting, we need $\beta$ to satisfy
\begin{align*}
c&\ge\max\left\{\frac{\delta}{L}, \frac{4\|L_{+}\|}{\sigma_{\min}(L_{-})}\right\}\\
\beta&> \frac{L}{2}\left(2c+1+\sqrt{(2c+1)^2+\frac{16 L^2}{\sigma_{\min}(L_{-})}}\right).
\end{align*}
Note that for EXTRA, a sufficient condition for  the penalty parameter is that $\beta> L$. Clearly the above requirement for $\beta$ is  at least $\frac{L}{2}\left(4+\sqrt{25+\frac{16 L^2}{\sigma_{\min}(L_{-})}}\right)$ larger. However, this is reasonable since Prox-GPDA is capable of dealing with nonconvex problems as well.

Note that our analysis is  by no means any extension to that of EXTRA. Fundamentally, the analysis of EXTRA relies upon showing the descent of the conventional potential function $\|x^*-x^r\|_{G}$ (i.e., the distance to the optimal solution set), where $G\succeq 0$ being some problem dependent matrix and $x^*$ is a solution in the {\it global} optimal solution set. In our nonconvex setting, such measure is not useful anymore. 

We remark that each node $i$ can distributedly implement iteration \eqref{eq:EXTRA} by performing the following
\begin{align}\label{eq:EXTRA:distributed}
x_i^{r+1}& = x_i^r -\frac{1}{2\beta d_i} \left(\nabla f_i(x_i^r) - \nabla f_i(x_i^{r-1})\right) + \sum_{j\in \cN(i)}\frac{1}{d_i} x_j^r -\frac{1}{2}\left(\sum_{j\in \cN(i)}\frac{1}{d_i} x_j^{r-1} + x_i^{r-1}\right)
\end{align}
where $\cN(i)$ denotes the set of neighbors for node $i$
\begin{align*}
\cN(i):=\{j\mid (i,j)\in \cE, i\ne j\}.
\end{align*}
Clearly, at iteration $r+1$, besides the local gradient information, node $i$ only needs the aggregated information from its neighbors, $\sum_{j\in \cN(i)}x^r_j$. Also such aggregated sum is required to be stored in local memory for at {\it least one more iteration}, because in order to carried out the $x^{r+1}_i$ update, node $i$ also needs $\sum_{j\in \cN(i)}x^{r-1}_j$.

We also remark that one can apply Prox-PDA (Algorithm 1) to the distributed setting as well. The resulting iteration is given below
\begin{align}\label{eq:EXTRA:alg1}
x^{r+1} & = x^r -\frac{1}{2\beta} D^{-1}\left(\nabla f(x^{r+1}) - \nabla f(x^{r})\right) + W x^r -\frac{1}{2}(I+W)x^{r-1}.
\end{align}
This algorithm can be implemented as follows. At iteration $1$, assuming that $x^0$ and $x^{-1}$ have been properly initialized. We generate $x^1$ according to the following
\begin{align}\label{eq:EXTRA:alg1:imp}
x^{1} & = x^0 -\frac{1}{2\beta} D^{-1}\left(\nabla f(x^{1}) - \nabla f(x^{0})\right) + W x^0-\frac{1}{2}(I+W)x^{-1}.
\end{align}
Equivalently, each node generates $x^1_i$ according to the following 
\begin{align}\label{eq:EXTRA:distributed:2}
x_i^{1}& = x_i^0 -\frac{1}{2\beta d_i} \left(\nabla f_i(x_i^{1}) - \nabla f_i(x_i^{0})\right) + \sum_{j\in \cN(i)}\frac{1}{d_i} x_j^0 -\frac{1}{2}\left(\sum_{j\in \cN(i)}\frac{1}{d_i} x_j^{-1} + x_i^{-1}\right)
\end{align}
To obtain such $x^{1}_i$, each node solves the following optimization problem
\begin{align}\label{eq:alg1:dist}
\min_{x_i} \; f_i(x_i)+{\beta d_i}\|x_i -\ell^1_i \|^2
\end{align}
where $\ell^1_i$ is given by
\begin{align}\label{eq:d}
\ell^1_i:= x^0_i+\frac{1}{2\beta d_i}\nabla f(x^0_i) +  \sum_{j\in \cN(i)}\frac{1}{d_i} x_j^0 -\frac{1}{2}\left(\sum_{j\in \cN(i)}\frac{1}{d_i} x_j^{-1} + x_i^{-1}\right).
\end{align}
Then at iteration $r+1$, node $i$ performs the same minimization as in \eqref{eq:alg1:dist}, only with $\ell^1_i$ replaced by $\ell^{r+1}_i$, which is similarly defined as in \eqref{eq:d}.

Also, we remark that applying either Algorithm 1 or Algorithm 2, the penalty parameter $\beta$ will still be selected according to \eqref{eq:c:bound} and \eqref{eq:beta:bound}. 
Intuitively, $\beta$ is decreasing with respect to the smallest nonzero eigenvalue of $L_{-}$, increasing with the maximum eigenvalue of $L_{+}$, and finally proportional to $L = \max_{i}L_i$.

Finally, in applications where explicitly selecting the stepsizes is difficult, an alternative is to use the Prox-GPDA-IP. To derive the iterates, let us select the $B$ matrix such that $B^T B = L_{+} +I_N$ (in order to satisfy \eqref{eq:B:inc}). Using this choice of $B$, the optimality condition for the $x$-subproblem \eqref{eq:x:update:inc:2} is given by
\begin{align*}
\nabla f(x^r) +A^T \mu^r +\beta^{r+1} (2 D + I_N) x^{r+1} - \beta^{r+1} (L_{+}+I_N)x^r=0.
\end{align*} 
Subtracting the same equation evaluated at the previous iteration, we obtain
\begin{align*}
&\nabla f(x^r) - \nabla f(x^{r-1}) + \beta^r L_{-} x^r +\beta^{r+1} (2 D +I_N) (x^{r+1}-x^r) \nonumber\\
&\quad\quad\quad\quad\quad -(\beta^r-\beta^{r+1})(2D+I_N)x^r-  (L_{+}+I_N)\left(\beta^{r+1}x^{r}-\beta^r x^{r-1}\right) = 0
\end{align*}
where we have used the fact that $A^T (\mu^{r}-\mu^{r-1})=\beta ^{r}A^T A x^r = \beta^r L_{-}x^r$. Rearranging terms, we have
\begin{align}\label{eq:EXTRA:increase}
x^{r+1} &= x^r -\frac{1}{2\beta^{r+1}} (D+\frac{1}{2}I_N)^{-1}\left(\nabla f(x^r) - \nabla f(x^{r-1})\right) + \frac{1}{2} (D+\frac{1}{2}I_N)^{-1} (L_{+}-L_{-}+I_N)x^r\nonumber\\
&\quad  -\frac{\beta^r-\beta^{r+1}}{2\beta^{r+1}}(D+\frac{1}{2}I_N)^{-1}L_{-}x^r-\frac{\beta^r}{2\beta^{r+1}} (D+\frac{1}{2}I_N)^{-1} (L_{+}+I_N)x^{r-1}-\frac{\beta^{r+1}-\beta^r}{\beta^{r+1}}x^r\nonumber\\
& = x^r -\frac{1}{2\beta^{r+1}} (D+\frac{1}{2}I_N)^{-1}\left(\nabla f(x^r) - \nabla f(x^{r-1})\right) + W x^r -\frac{\beta^{r+1}-\beta^r}{\beta^{r+1}}x^r\nonumber\\ &\quad -\frac{\beta^r-\beta^{r+1}}{2\beta^{r+1}}(I-W)x^r-\frac{\beta^r}{2\beta^{r+1}}(I+W)x^{r-1} 
\end{align}
where in the last equality we have defined the {\it weight matrix} $W: = \frac{1}{2}(D+\frac{1}{2}I_N)^{-1} (L_{+}-L_{-}+I_N)$, which is a row stochastic matrix.


The discussion in this section is summarized in the following corollary . 
\begin{corollary}\label{cor:distribtued}
	{\it Consider the distributed optimization problem \eqref{eq:global:consensus:equiv}. Suppose that the graph $(\cV, \cE)$ is connected. Then we have the following claims. 
		\begin{enumerate}
			\item Suppose Assumption A is satisfied. Suppose $\beta$ and $c$ satisfy \eqref{eq:beta:bound} and \eqref{eq:c:bound}. Then all the conclusions in Theorem  \ref{thm:main} hold true for the distributed iterations \eqref{eq:EXTRA} and \eqref{eq:EXTRA:alg1}, respectively.
			
			\item Suppose Assumption A is satisfied, and $\{\beta^r\}$ satisfies \eqref{eq:beta}. Then all the conclusions in Theorem \ref{thm:inc} hold true for the distributed iteration \eqref{eq:EXTRA:increase}. 
			 		\item Moreover, in either case, a stationary solution of the original unconstrained problem \eqref{eq:global:consensus} is achieved. 
		\end{enumerate}	}

\end{corollary}
The first and the second statements  directly follow the results in Theorem \ref{thm:main}, Corollary \ref{cor:proximal} and Theorem \ref{thm:inc}. The third statement is easy to see because the stationary solution of problem \eqref{eq:global:consensus:equiv} is given by
\begin{align*}
\nabla  f(x^*)+A^T \mu^* = 0, \quad A x^* = 0.
\end{align*}
Left multiplying the first equality by an all one vector $1\in\mathbb{R}^{N}$, we have
\begin{align*}
\sum_{i=1}^{N}\nabla f_i(x^*_i) + (A \times 1)^T \mu^* = \sum_{i=1}^{N}\nabla f_i(x^*_i) =0,
\end{align*}
where the first equality is true because $1$ is in the null space of the incidence  matrix $A$. The fact that $A x^* = 0$ implies that $x_i=x_j$, for all $i\ne j$. Therefore we conclude that
\begin{align*}
\sum_{i=1}^{N}\nabla f_i(x^*_j) =0,\quad j=1,\cdots, N.
\end{align*}

This corollary suggests that iterations \eqref{eq:EXTRA} and \eqref{eq:EXTRA:alg1} also achieve a global sublinear convergence. Note that there has been a few recent works on distributed nonconvex optimization, for example \cite{Lorenzo16,Tatarenko05} and the references therein. These works design algorithms under different assumptions on the problem as well as on the network structure. However, central to these works is the use of certain diminishing stepsize for th local updates, which results in no global convergence rate guarantees. To the best of our knowledge, the Prox-PDA based distributed algorithms are the first ones that provably achieve global sublinear convergence rate for nonconvex distributed optimization. 

\section{Generalization: Distributed Nonconvex Matrix Factorization}
In this section we study a variant of the Prox-PDA algorithm for a distributed matrix factorization problem. 

Consider the following matrix factorization problem \cite{Ling12matrix}
\begin{align}\label{eq:mat:problem}
\min_{X, Y}\; \frac{1}{2}\|X Y -Z\|_F^2 + \gamma \|X\|^2_F+ h(Y)=\sum_{i=1}^{N}\frac{1}{2}\|X y_i -z_i\|^2+\gamma \|X\|^2_F+h_i(y_i), \; \st \;  \|y_i\|^2\le \tau
\end{align}
where $X\in\mathbb{R}^{M\times K}$, $Y\in\mathbb{R}^{K\times N}$, and $y_i\in\mathbb{R}^{K\times 1}$ is the $i$th column of $Y$, respectively; $Z\in\mathbb{R}^{M\times N}$ is the observation matrix; $h(Y):=\sum_{i=1}^{N}h_i(y_i)$ is some convex but possibly nonsmooth penalization term; $\gamma>0$ is some given constant. We assume that $h(Y)$ is lower bounded over $\mbox{dom}\;(h)$. One application of problem \eqref{eq:mat:problem} is the distributed {\it sparse dictionary learning} problem where $X$ is the dictionary to be learned, each $z_i$ is a training data sample, and each $y_i$ is the sparse coefficient corresponding to the particular training sample $z_i$.  The constraint $\|y_i\|^2$ simply says that the size of the coefficient must be bounded. 

Consider a distributed scenario where $N$ agents form a graph $\{\cV,\cE\}$, each having a column of $Y$ (note, this can be easily generalized to the case where a subset of columns are available for each agent), we reformulate problem \eqref{eq:mat:problem} as
\begin{align*}
\min_{\{X_i\},\{y_i\}}&\; \sum_{i=1}^{N}\left(\frac{1}{2}\|X_i y_i - z_i\|^2+h_i(y_i) + \gamma \|X_i\|^2_F\right)\nonumber\\
\st&\; \|y_i\|^2\le \tau, \quad X_i = X_j\; \forall~(i,j)\in \cE.
\end{align*}

Stacking all the variables $X_i$, let us define $\bX := [X_1; X_2; \cdots, X_N]\in\mathbb{R}^{NM\times K}$. Define the block signed incidence matrix as $\bA = A\otimes I_M\in\mathbb{R}^{E M\times NM}$, where $A$ is the standard graph incidence matrix. Define the block  signless incidence matrix $\bB\in\mathbb{R}^{EM\times NM}$ similarly.
If the graph is connected, then the following condition implies network-wide consensus
$$\bA\bX =\bzero.$$
Using this condition, we  formulate the distributed matrix factorization problem as 
\begin{align}\label{eq:dis:matrix}
\begin{split}
\min_{\{X_i\},\{y_i\}}&\; f(\bX, Y) + h(Y):=\sum_{i=1}^{N}\frac{1}{2}\|X_i y_i - z_i\|^2+\gamma \|X_i\|^2_F+ h_i(y_i)\\
\st&\; \|y_i\|^2\le \tau, \quad \bA \bX =\bzero.
\end{split}
\end{align}
Clearly the above problem does not fall into the form of \eqref{eq:original}, because there are two block variables $\{X_i\}$ and $\{y_i\}$ in the objective, but the linear constraint only has to do with the $X$-block.  Moreover, the objective function couples among the variable blocks $\{X_i\}$ and $\{y_i\}$ in a nonconvex manner, and neither $\{X_i\}$ nor $\{y_i\}$  has Lipschitz continuous gradient. The latter fact poses significant difficulty in algorithm development and analysis.

Define the block-signed and the block-signless Laplacians as follows
$$\bL_{-} = \bA^T\bA\in\mathbb{R}^{NM\times NM}, \quad \bL_{+} = \bB^T \bB\in\mathbb{R}^{NM\times NM}.$$

The augmented Lagrangian for the above problem is given by
\begin{align}\label{eq:mat:lag}
L_{\beta}(\bX, Y, \bfOmega) = \sum_{i=1}^{N}\left(\frac{1}{2}\|X_i y_i - z_i\|^2+\gamma\|X_i\|^2_F+h_i(y_i)\right)+\langle \bfOmega, \bA\bX\rangle  +\frac{\beta}{2}\langle \bA\bX, \bA\bX\rangle,
\end{align}
where $\bfOmega:=\{\Omega_e\}\in\mathbb{R}^{EM\times K}$ is the dual variable, where for each $e=(i,j)$,  $\Omega_e\in\mathbb{R}^{M\times K}$ represents the dual variable for the consensus constraint $X_i = X_j$.

Let us consider the following generalization of Algorithm 1 for distributed matrix factorization.

\begin{center}
	\fbox{
		\begin{minipage}{5in}
			\smallskip
			\centerline{\bf {Algorithm 6. Prox-PDA for Distributed Matrix Factorization}}
			\smallskip
			At iteration $0$, initialize $\bfOmega^0=\mathbf{0}$, and $\bX^0, y^0$;
			
			At each iteration $r+1$, update variables  by:
			\begin{subequations}
				\begin{align}
				\theta^r_i& = \|X^r_i y^r_i - z_i\|^2,\quad\forall~i;\label{eq:theta:update}\\
				y_i^{r+1} & =\arg\min_{\|y_i\|^2\le \tau}\frac{1}{2}\|X^r_i y_i-z_i\|^2+h_i(y_i)+\frac{\theta^r_i}{2}\|y_i-y^r_i\|^2,\quad \forall~i;\label{eq:y:update:matrix}\\
				\bX^{r+1}& =\arg\min_{\bX\in\mathbb{R}^{NM\times K}}\; f(\bX, Y^{r+1})+\langle\bfOmega^r, \bA\bX \rangle;\nonumber\\
				&\quad +\frac{\beta}{2}\langle \bA\bX, \bA\bX\rangle +\frac{\beta}{2}\langle \bB(\bX-\bX^r), \bB (\bX-\bX^r)\rangle;\label{eq:x:update:matrix}\\
				\bfOmega^{r+1}& = \bfOmega^r +\beta \bA \bX^{r+1}.				\label{eq:omega:update:matrix}
				\end{align}
			\end{subequations}
			
		\end{minipage}
	}
\end{center}
\vspace{0.2cm}
In the above algorithm we have introduced a new sequence $\theta^r_i\ge 0$, which is some iteration-dependent coefficient representing the size of the local factorization error. We note that including the proximal term $\frac{\theta^r_i}{2}\|y_i-y^r_i\|^2$ is the key to achieving convergence for Algorithm 6. Ideally, the effect of such proximal term will disappear as the algorithm approaches convergence, since we would expect that the local factorization error becomes small. Again one should note that $\frac{\beta}{2}\langle \bA\bX, \bA\bX\rangle +\frac{\beta}{2}\langle \bB(\bX-\bX^r), \bB (\bX-\bX^r)\rangle$ is strongly convex in $\bX$. 

We briefly comment on how the algorithm can be implemented in a distributed manner. First note that the $y$ subproblem \eqref{eq:y:update:matrix} is naturally distributed to each node, that is, only local information is needed to perform the update. Second, the $\bX$ subproblem \eqref{eq:x:update:matrix} can also be decomposed into $N$  subproblems, one for each node. To be more specific, let us examine the terms in \eqref{eq:x:update:matrix} one by one. First, the term $f(\bX, Y^{r+1})=\sum_{i=1}^{N}\left(\frac{1}{2}\|X_i y^{r+1}_i-z_i\|^2+h_i(y_i)+\gamma\|X_i\|^2_F\right)$, hence it is decomposable. Second, the term $\langle\bfOmega^r, \bA\bX \rangle$ can be expressed as
\begin{align*}
\langle\bfOmega^r, \bA\bX \rangle = \sum_{i=1}^{N}\sum_{e\in U(i)}\langle  \Omega^r_e  ,X_i \rangle - \sum_{e\in H(i)}\langle  \Omega^r_e  ,X_i \rangle 
\end{align*}
where the sets $U(i)$ and $H(i)$ are defined as
\begin{align*}
U(i):= \{e \mid e=(i,j)\in\cE, i\ge j\},\quad H(i):= \{e \mid e=(i,j)\in\cE, j\ge i\}.
\end{align*}
Similarly, we have
\begin{align}
\langle \bB\bX^r, \bB \bX \rangle & =\sum_{i=1}^{N} \left\langle X_i, d_i X^r_i+\sum_{j\in N(i)} X^r_j\right\rangle \nonumber\\
\frac{\beta}{2}\left(\langle \bA \bX, \bA \bX\rangle+\langle \bB \bX, \bB \bX\rangle\right) & = \beta \langle \bD \bX, \bX\rangle  =\beta\sum_{i=1}^{N}d_i \|X_i\|^2_F  \nonumber
\end{align}
where $\bD = D\otimes I_{M}\in\mathbb{R}^{NM\times NM}$, where $D=\mbox{diag}[d_1,\cdots, d_n]\in\mathbb{R}^{N\times N}$ is the degree matrix. 
Therefore it is easy to see that the $\bX$ subproblem \eqref{eq:x:update:matrix} is separable over the distributed agents. 

Finally, one can verify that the $\bfOmega$ update step \eqref{eq:omega:update:matrix} can be implemented by each edge $e\in \cE$ as follows
\begin{align*}
\Omega^{r+1}_e = \Omega^{r}_e + \beta\left(X^{r+1}_i-X^{r+1}_j\right), \quad e=(i,j), i\ge j. 
\end{align*}

To show convergence rate of the algorithm, we need the following definition. 
\begin{align*}
Q(\bX^{r+1}, Y^{r+1},\bfOmega^{r}):=\beta\|\bA \bX^{r+1}\|^2+ \|[\bZ^{r+1}_1; \bZ^{r+1}_2]\|^2
\end{align*}
where 
\begin{align}
\bZ^{r+1}_1 &: = \nabla_{\bX}L(\bX^{r+1}, Y^{r+1},\bfOmega^{r}); \nonumber\\
\bZ^{r+1}_2 &: = Y^{r+1}-\prox_{h+\iota(\mathcal{Y})}\left[Y^{r+1}-\nabla_{Y}\left(L(\bX^{r+1}, Y^{r+1},\bfOmega^{r}) -h(Y)\right)\right].\nonumber
\end{align}
In the above derivation, we have defined the proximity operator for a given convex lower semi-continuous function $p(\cdot)$ as 
\begin{align}\label{eq:prox}
\prox_{h}(c) = \arg\min_{z}\;p(z)+\frac{1}{2}\|z-c\|^2.
\end{align}
We have used $\mathcal{Y}:=\bigcup_{i}\left\{\|y_i\|^2\le \tau\right\}$  to denote the feasible set of $Y$, and used $\iota(\mathcal{Y})$ to denote the indicator function of such set.
Similarly as in Section \ref{sec:analysis}, we can show that $Q(\bX^{r+1}, Y^{r+1},\bfOmega^{r+1})\to 0$ implies that every limit point of $(\bX^{r+1},Y^{r+1}, \bfOmega^{r+1})$ is a stationary point of problem \eqref{eq:dis:matrix}. 

Next we present the main convergence analysis for Algorithm 6. The proof is long therefore we delegate it to Section \ref{sec:proof:matrix}.  

\begin{theorem}\label{thm:matrix}
Consider using Algorithm 6 to solve the distributed matrix factorization problem \eqref{eq:dis:matrix}. Suppose that $h(Y)$ is lower bounded over $\mbox{dom}\; h(x)$, and that the penalty parameter $\beta$, together with two positive constants $c$ and $d$,  satisfies the following conditions
\begin{align}\label{eq:condition:matrix:thm}
\begin{split}
\frac{\beta+2\gamma}{2}-\frac{8(\tau^2+4\gamma^2)}{\beta{\sigma_{\min}(A^T A)}}-\frac{cd}{2}&>0, \quad \frac{1}{2}-\frac{{8}}{\sigma_{\min}(A^T A)\beta}-\frac{c}{d}>0\\
\frac{1}{2}-\frac{{8\tau}}{\sigma_{\min}(A^T A)\beta}-\frac{c\tau}{d}&>0 ,\quad \frac{c\beta}{2}-\frac{2\beta\|B^T B\|}{{\sigma_{\min}(A^T A)}}>0.
\end{split}
\end{align}
Then in the limit, consensus will be achieved, i.e., 
$$X_i = X_j, \quad \forall~(i,j)\in \cE.$$
Further, the sequences $\{\bX^{r+1}\}$ and $\{\bfOmega^{r+1}\}$ are both bounded, and every limit point generated by Algorithm 6 converges to a stationary point for problem \eqref{eq:mat:problem}. 

Additionally, Algorithm 6 converges sublinearly, i.e., the measure  $Q(\bX^{r+1},Y^{r+1},\bfOmega^{r+1})$ decreases to $0$ in the same manner as in Theorem \ref{thm:main}. Specifically,  for any given $\varphi>0$, define $T$ to be the first time that the optimality gap reaches below $\varphi$, i.e., 
$$T:=\arg\min_{r}\; Q(\bX^{r+1}, Y^{r+1}, \bfOmega^{r})\le \varphi.$$
Then there exists a constant $\nu>0$ such that the following is true
\begin{align*}
\varphi\le \frac{\nu}{T-1}.
\end{align*}

\end{theorem}

We can see that it is always possible to find the tuple $\{\beta, c, d>0\}$ that satisfies \eqref{eq:condition:matrix:thm}. For example $c$ can be solely determined by the last inequality; $d$ needs to be chosen large enough such that $1/2-\frac{c}{d}>0$ and $1/2-\frac{c\tau}{d}>0$. After $c$ and $d$ are fixed, one can always choose $\beta$ large enough to satisfy the first three conditions. In practice, we typically prefer to choose $\beta$ as small as possible to improve the convergence speed. Therefore empirically one can start with (for some small $\nu>0$)
$$c= \frac{4\|B^T B\|}{\sigma_{\min}(A^T A)}+\nu \quad d=\max\{4, 2 c \tau \}$$
and then gradually increase $d$ to find an appropriate $\beta$ that satisfies the first three conditions. Of course, one also has the option of utilizing increasing penalty parameters, just as what we have done in Section \ref{sec:proof:increasing}. 

\section{Concluding Remarks}
In this paper, we have proposed a decomposition approach for certain  linearly constrained nonconvex smooth optimization problem. Our developed algorithms, mainly based upon a novel proximal primal-dual augmented Lagrangian method, are  able to decompose the optimization variables, resulting in simple subproblems that can often be solved in closed-form. By constructing a new potential function, which is a conic combination of the augmented Lagrangian, the size of the constraint violation and q certain proximal term, we have shown that the proposed Prox-PDA and its various extensions converge globally sublinearly to the set of stationary solutions. Surprisingly, when specializing a variant of Prox-PDA to the nonconvex distributed optimization problem, the proposed algorithm recovers the popular EXTRA algorithm, indicating that such algorithm converges globally sublinearly even for nonconvex problems. 

The proposed proximal primal-dual based algorithm can have many extensions or generalizations. In the paper we have discussed one such generalization to a (distributed) matrix factorization problem. Can we deal with nonsmooth nonconvex terms in the objective, such as indicator functions of convex/nonconvex sets? Can we randomized the algorithm so that each time a randomly selected subproblem is solved instead of the full subproblems? Can we apply our approach to stochastic optimization problems where the objective function involves the expectation of certain nonconvex function? Can we show that some variant of the Prox-PDA converges to local optimal solutions instead of stationary solutions? These are all very interesting research questions that require further investigation. 

\section{Proof of Convergence for Algorithm 4}\label{sec:proof:increasing}
Our analysis consists of a series of steps. 

\noindent{\bf Step 1.} Our first step is again to bound the size of the successive difference of $\{\mu^r\}$. To this end, write down the optimality condition for the $x$-update \eqref{eq:x:update:inc} as
\begin{align}\label{eq:mu:update:inc2}
A^T \mu^{r+1} = -\nabla f(x^{r+1}) -\beta^{r+1} B^T B (x^{r+1}-x^r).
\end{align}
Subtracting the previous iteration, we obtain
\begin{align}\label{eq:differ:mu:inc}
A^T (\mu^{r+1} -\mu^r) &= -(\nabla f(x^{r+1}) - \nabla f(x^{r})) -\beta^{r} B^T B \left((x^{r+1}-x^r)-(x^{r}-x^{r-1})\right) \nonumber\\
&\quad - (\beta^{r+1}-\beta^r)B^T B(x^{r+1}-x^r).
\end{align}
Therefore, we have
\begin{align}
\frac{1}{\beta^{r+1}}\|\mu^{r+1}-\mu^r\|^2 &\le\frac{3}{\beta^{r+1}\sigma_{\min}(A^T A)}\left(L^2+(\beta^{r+1}-\beta^r)^2\|B^T B\|\right)\|x^{r+1}-x^r\|^2\nonumber\\
&\quad + \frac{3 (\beta^r)^2}{\beta^{r+1}\sigma_{\min}(A^T A)}\left\|B^T B \left((x^{r+1}-x^r)-(x^{r}-x^{r-1})\right)\right\|^2\label{eq:mu:bound:inc}.
\end{align}

Also from the optimality condition we have the following relation
\begin{align}
x^{r+1} = x^r - \frac{1}{\beta^{r+1}}(B^T B)^{-1}\left(\nabla f(x^{r+1})+A^T \mu^{r+1}\right): =  x^r - \frac{1}{\beta^{r+1}}v^{r+1}\label{eq:def:v}.
\end{align}
where we have defined the primal update direction $v^{r+1}$ as 
\begin{align*}
v^{r+1} = (B^T B)^{-1}\left(\nabla f(x^{r+1})+A^T \mu^{r+1}\right).
\end{align*}

\noindent{\bf Step 2.} In the second step we analyze the descent of the augmented Lagrangian. We have the following estimate 
\begin{align}\label{eq:lag:inc}
&L_{\beta^{r+1}}(x^{r+1},\mu^{r+1}) - L_{\beta^r}(x^r,\mu^r)\nonumber\\
& = L_{\beta^{r+1}}(x^{r+1},\mu^{r+1}) - L_{\beta^{r+1}}(x^{r+1},\mu^r) + L_{\beta^{r+1}}(x^{r+1},\mu^{r}) - L_{\beta^{r+1}}(x^r,\mu^r)+ L_{\beta^{r+1}}(x^{r},\mu^{r}) - L_{\beta^r}(x^r,\mu^r)\nonumber\\
&\stackrel{\rm(i)}\le \frac{1}{\beta^{r+1}}\|\mu^{r+1}-\mu^r\|^2 + \frac{\beta^{r+1}-\beta^r}{2(\beta^r)^2}\|\mu^{r}-\mu^{r-1}\|^2-\frac{\beta^{r+1}-L}{2}\|x^{r+1}-x^r\|^2\nonumber\\
&\stackrel{\rm(ii)}\le -\left(\frac{\beta^{r+1}-L}{2} - \frac{3}{\beta^{r+1}\sigma_{\min}(A^T A)}\left(L^2+(\beta^{r+1}-\beta^r)^2\|B^T B\|\right)\right)\|x^{r+1}-x^r\|^2 +  \frac{\beta^{r+1}-\beta^r}{2(\beta^r)^2}\|\mu^{r}-\mu^{r-1}\|^2\nonumber\\
&\quad +\frac{3 (\beta^r)^2}{\beta^{r+1}\sigma_{\min}(A^T A)}\left\|B^T B \left((x^{r+1}-x^r)-(x^{r}-x^{r-1})\right)\right\|^2
\end{align}
where in ${\rm (i)}$ we have used the optimality of the $x$-subproblem (cf. the derivation in \eqref{eq:lag:bound:derive}); in ${\rm (ii)}$ we have applied \eqref{eq:mu:bound:inc}. 

\noindent{\bf Step 3.} In the third step, we construct the remaining part of the potential function. We have the following two inequalities from the optimality condition of the $x$-update \eqref{eq:x:update:inc}
\begin{align}
\left\langle  \nabla f(x^{r+1}) + A^T \mu^{r+1}+\beta^{r+1}B^T B(x^{r+1}-x^r), x^{r+1}-x  \right\rangle\le 0, \; \forall ~x\in\mathbb{R}^N\nonumber\\
\left\langle  \nabla f(x^{r}) + A^T \mu^{r}+\beta^{r}B^T B(x^{r}-x^{r-1}), x^{r}-x  \right\rangle\le 0, \; \forall ~x\in\mathbb{R}^N\nonumber.
\end{align}
Plugging $x=x^r$ and $x=x^{r+1}$ to these two equations and adding them together, we obtain\label{eq:second:estimate:inc}
\begin{align}
&\langle A^T (\mu^{r+1}-\mu^r), x^{r+1}-x^r\rangle\nonumber\\
& \le -\langle \nabla f(x^{r+1})-\nabla f(x^r), x^{r+1}-x^r\rangle -\langle B^T B(\beta^{r+1}(x^{r+1}-x^r) -\beta^r(x^{r}-x^{r-1})), x^{r+1}-x^r \rangle. \nonumber
\end{align}
The lhs of the above inequality can be expressed as 
\begin{align}
&\langle A^T (\mu^{r+1}-\mu^r), x^{r+1}-x^r\rangle\nonumber\\
& =\frac{\beta^{r+1}}{2}\left(\|Ax^{r+1}-b\|^2- \|Ax^{r}-b\|^2+ \|A (x^{r+1}-x^r)\|^2\right)\nonumber\\
& =\frac{\beta^{r+1}}{2}\|Ax^{r+1}-b\|^2- \frac{\beta^r}{2}\|Ax^{r}-b\|^2+ \frac{\beta^{r+1}}{2}\|A (x^{r+1}-x^r)\|^2+\frac{\beta^{r}-\beta^{r+1}}{2}\|Ax^r-b\|^2\nonumber.
\end{align}
The rhs of \eqref{eq:second:estimate:inc} can be bounded as
\begin{align}
&  -\langle \nabla f(x^{r+1})-\nabla f(x^r), x^{r+1}-x^r\rangle -\langle B^T B(\beta^{r+1}(x^{r+1}-x^r) -\beta^r(x^{r}-x^{r-1})), x^{r+1}-x^r \rangle\nonumber\\
&\le L\|x^{r+1}-x^r\|^2-(\beta^{r+1}-\beta^r)\|x^{r+1}-x^r\|^2_{B^T B}\nonumber\\
&\quad +\frac{\beta^r}{2}\left(\|x^{r}-x^{r-1}\|^2_{B^T B} -\|x^{r}-x^{r+1}\|^2_{B^T B}-\|(x^r-x^{r-1})-(x^{r+1}-x^r)\|^2_{B^T B}\right) \nonumber\\
&= L\|x^{r+1}-x^r\|^2-\frac{\beta^{r+1}-\beta^r}{2}\|x^{r+1}-x^r\|^2_{B^T B}\nonumber\\
&\quad +\frac{\beta^r}{2}\|x^{r}-x^{r-1}\|^2_{B^T B} -\frac{\beta^{r+1}}{2}\|x^{r}-x^{r+1}\|^2_{B^T B}-\frac{\beta^r}{2}\|(x^r-x^{r-1})-(x^{r+1}-x^r)\|^2_{B^T B} \nonumber\\
&\stackrel{\eqref{eq:beta}}\le L\|x^{r+1}-x^r\|^2+\frac{\beta^r}{2}\|x^{r}-x^{r-1}\|^2_{B^T B} -\frac{\beta^{r+1}}{2}\|x^{r}-x^{r+1}\|^2_{B^T B}-\frac{\beta^r}{2}\|(x^r-x^{r-1})-(x^{r+1}-x^r)\|^2_{B^T B} \nonumber.
\end{align}
Therefore, combining the above three inequalities we obtain
\begin{align}
&\frac{\beta^{r+1}}{2}\|Ax^{r+1}-b\|^2 + \frac{\beta^{r+1}}{2}\|x^{r}-x^{r+1}\|^2_{B^T B}\nonumber\\
&\le \frac{\beta^r}{2}\|Ax^{r}-b\|^2 +\frac{\beta^r}{2}\|x^{r}-x^{r-1}\|^2_{B^T B} +\frac{\beta^{r+1}-\beta^{r}}{2(\beta^r)^2}\|\mu^{r-1}-\mu^r\|^2+  L\|x^{r+1}-x^r\|^2\nonumber\\
&\quad-\frac{\beta^r}{2}\|(x^r-x^{r-1})-(x^{r+1}-x^r)\|^2_{B^T B}\nonumber .
\end{align}
Multiplying both sides by $\beta^{r}$, we obtain
\begin{align}\label{eq:bound:key:inc}
&\frac{\beta^{r+1}\beta^r}{2}\|Ax^{r+1}-b\|^2 + \frac{\beta^{r+1}\beta^r}{2}\|x^{r}-x^{r+1}\|^2_{B^T B}\nonumber\\
&\le \frac{\beta^r\beta^{r-1}}{2}\|Ax^{r}-b\|^2 +\frac{\beta^r\beta^{r-1}}{2}\|x^{r}-x^{r-1}\|^2_{B^T B} +\frac{\beta^{r+1}-\beta^{r}}{2\beta^r}\|\mu^{r-1}-\mu^r\|^2+  \beta^r L\|x^{r+1}-x^r\|^2\nonumber\\
&\quad-\frac{(\beta^r)^2}{2}\|(x^r-x^{r-1})-(x^{r+1}-x^r)\|^2_{B^T B} \nonumber\\
&\quad +\frac{\beta^r(\beta^r-\beta^{r-1})}{2}\|Ax^{r}-b\|^2 +\frac{\beta^r(\beta^r-\beta^{r-1})}{2}\|x^{r}-x^{r-1}\|_{B^T B}^2 \nonumber\\
&= \frac{\beta^r\beta^{r-1}}{2}\|Ax^{r}-b\|^2 +\frac{\beta^r\beta^{r-1}}{2}\|x^{r}-x^{r-1}\|^2_{B^T B} +\frac{\beta^{r+1}-\beta^{r-1}}{2\beta^r}\|\mu^{r-1}-\mu^r\|^2+  \beta^r L\|x^{r+1}-x^r\|^2\nonumber\\
&\quad -\frac{(\beta^r)^2}{2}\|(x^r-x^{r-1})-(x^{r+1}-x^r)\|^2_{B^T B} +\frac{\beta^r(\beta^r-\beta^{r-1})}{2}\|x^{r}-x^{r-1}\|_{B^T B}^2 
\end{align}
where in the last equality we have merged the terms $\frac{\beta^{r+1}-\beta^{r}}{2\beta^r}\|\mu^{r-1}-\mu^r\|^2$ and $\frac{\beta^r(\beta^r-\beta^{r-1})}{2}\|Ax^{r}-b\|^2$.

\noindent{\bf Step 4.} In this step we construct and estimate the descent of the potential function. For some given $c>0$, let us define the potential function as
\begin{align}
P_{\beta^{r+1}, c}(x^{r+1}, x^r, \mu^{r+1})= L_{\beta^{r+1}}(x^{r+1},\mu^{r+1}) + \frac{c\beta^{r+1}\beta^r}{2}\|Ax^{r+1}-b\|^2 + \frac{c\beta^{r+1}\beta^r}{2}\|x^{r}-x^{r+1}\|^2_{B^T B}.\nonumber
\end{align}
Note that this potential function has some major differences compared with the one we used before; cf. \eqref{eq:potential}. In particular, the second and the third terms are now quadratic, rather than linear,  in the penalty parameters. This new construction is the key to our following analysis. 

Then combining the estimate in \eqref{eq:bound:key:inc} and \eqref{eq:lag:inc}, we obtain
\begin{align}\label{eq:potential:descent:inc}
&P_{\beta^{r+1}, c}(x^{r+1}, x^r, \mu^{r+1}) - P_{\beta^{r}, c}(x^{r}, x^{r-1}, \mu^{r})\nonumber\\
&\le -\left(\frac{\beta^{r+1}-L}{2} - \frac{3}{\beta^{r+1}\sigma_{\min}(A^T A)}\left(L^2+(\beta^{r+1}-\beta^r)^2\|B^T B\|\right)-c\beta^r L\right)\|x^{r+1}-x^r\|^2 \nonumber\\
&\quad +  \frac{\beta^{r+1}-\beta^{r-1}}{2\beta^r}(\frac{1}{\beta^r}+c)\|\mu^{r}-\mu^{r-1}\|^2+\frac{c\beta^r(\beta^r-\beta^{r-1})}{2}\|x^{r}-x^{r-1}\|_{B^TB}^2 \nonumber\\
&\quad -\left(\frac{c(\beta^{r})^2}{2}- \frac{3 (\beta^r)^2\|B^T B\|}{\beta^{r+1}\sigma_{\min}(A^T A)}\right)\left\|(x^{r+1}-x^r)-(x^{r}-x^{r-1})\right\|_{B^T B}^2
\end{align}
where in the inequality we have also used the fact that $\beta^{r}\ge \beta^{r-1}$. 

Taking the sum of $r$ from $t$ to $T+1$ (for some $T>t>1$) and utilize again the estimate in \eqref{eq:mu:bound:inc}, we have 
	\begin{align}\label{eq:potential:sum:descent}
	&P_{\beta^{T+1}, c}(x^{T+1}, x^T, \mu^{T+1}) - P_{\beta^{t}, c}(x^{t}, x^{t-1}, \mu^{t})\nonumber\\
	&\le -\sum_{r=t}^{T}\bigg(\frac{\beta^{r+1}-L}{2} - \frac{3+3(1/\beta^r+c)(\beta^{r+1}-\beta^{r-1})/2\beta^r}{\beta^{r+1}\sigma_{\min}(A^T A)}\left(L^2+(\beta^{r+1}\nonumber-\beta^{r-1})^2\|B^T B\|\right)\nonumber\\
	&\quad -c\beta^{r} L- \frac{c \beta^{r+1}(\beta^{r+1}-\beta^{r})\|B^T B\|}{2}\bigg)\|x^{r+1}-x^r\|^2 \nonumber\\
	&\quad -\left(\frac{c(\beta^{r+1})^2}{2}- \frac{(3+3(1/\beta^r+c)(\beta^{r+1}-\beta^{r-1})/2\beta^r) (\beta^r)^2\|B^T B\|}{\beta^{r+1}\sigma_{\min}(A^T A)}\right)\left\|(x^{r+1}-x^r)-(x^{r}-x^{r-1})\right\|_{B^T B}^2\nonumber\\
	& \quad  +\frac{c\beta^t(\beta^t-\beta^{t-1})}{2}\|x^{t}-x^{t-1}\|_{B^T B}^2+ \frac{\beta^{t+1}-\beta^{t-1}}{2\beta^t}(1/\beta^t+c)\|\mu^{t}-\mu^{t-1}\|^2.
	\end{align}
First, note that for any $c\in(0,1)$, the coefficient in front of $\left\|(x^{r+1}-x^r)-(x^{r}-x^{r-1})\right\|_{B^T B}^2$ becomes negative for sufficiently large (but finite) $t$. This is because $\{\beta^r\}\to \infty$, and that the first term in the parenthesis scales in $\mathcal{O}((\beta^r)^2)$ while the second term scales in $\mathcal{O}(1)$ .
For the first term to be negative, we need $c>0$ to be {\it small enough}  such that the following is true for large enough $r$
\begin{align}
\frac{\beta^{r+1}-L}{2} -c\beta^{r} L- \frac{c\beta^{r+1}(\beta^{r+1}-\beta^{r})\|B^T B\|}{2}>\frac{\beta^{r+1}}{24}. \nonumber
\end{align}
Suppose that $r$ is large enough so that $(\beta^{r+1}-L)/2> {\beta^{r+1}}/{3}$, or equivalently $\beta^{r+1}>3L$.
Also choose $c=\min\{1/(4L), 1/(12\omega \|B^T B\|)\}$, where $\omega$ is given in \eqref{eq:beta}. Then we have
\begin{align}\label{eq:beta:c:inc}
&\frac{\beta^{r+1}-L}{2} -c\beta^{r} L- \frac{c \beta^{r+1}(\beta^{r+1}-\beta^{r})\|B^TB\|}{2}>\frac{\beta^{r+1}}{3} -\frac{\beta^{r+1}}{4}- \frac{\beta^{r+1}}{24}=\frac{\beta^{r+1}}{24}.
\end{align}
For this given $c$, we can also show that the following is true for sufficiently large $r$
\begin{align}
&\frac{3+3(1/\beta^r+c)(\beta^{r+1}-\beta^{r-1})/2\beta^r}{\beta^{r+1}\sigma_{\min}(A^T A)}\left(L^2+(\beta^{r+1}-\beta^r)^2\|B^T B\|\right)\le \frac{\beta^{r+1}}{48}\nonumber\\
&\left(\frac{c(\beta^{r+1})^2}{2}- \frac{(3+3(1/\beta^r+c)(\beta^{r+1}-\beta^{r-1})/2\beta^r) (\beta^r)^2\|B^T B\|}{\beta^{r+1}\sigma_{\min}(A^T A)}\right)\ge \frac{(\beta^{r+1})^2}{48}\nonumber.
\end{align}
In conclusion we have that for sufficiently large but finite $t_0$, we have
\begin{align}\label{eq:clean}
&P_{\beta^{T+1}, c}(x^{T+1}, x^T, \mu^{T+1}) - P_{\beta^{t_0}, c}(x^{t_0-1}, x^{t_0}, \mu^{t_0})\nonumber\\
&\le -\sum_{r=t_0}^{T}\frac{\beta^{r+1}}{48}\|x^{r+1}-x^r\|^2 -\frac{(\beta^{r+1})^2}{48}\left\|(x^{r+1}-x^r)-(x^{r}-x^{r-1})\right\|_{B^T B}^2\nonumber\\
& \quad  +\frac{c\beta^{t_0}(\beta^{t_0}-\beta^{t_0-1})}{2}\|x^{t_0}-x^{t_0-1}\|_{B^T B}^2+ \frac{\beta^{t_0+1}-\beta^{t_0-1}}{2\beta^{t_0}}(1/\beta^{t_0}+c)\|\mu^{t_0}-\mu^{t_0-1}\|^2.
\end{align}

Therefore we conclude that if $\{\beta^{r+1}\}$ satisfies \eqref{eq:beta}, and for $c$ sufficiently small, there exits a finite $t_0>0$ such that for all $T>t_0$, the first two terms of the rhs of \eqref{eq:potential:sum:descent} is negative.

\noindent{\bf Step 5.} Next we show that the potential function must be lower bounded. Observe that the augmented Lagrangian is given by
\begin{align}
&L_{\beta^{r+1}}(x^{r+1}, \mu^{r+1})\nonumber\\
&=f(x^{r+1}) +\langle\mu^{r+1}, Ax^{r+1} - b\rangle+\frac{\beta^{r+1}}{2}\|Ax^{r+1}-b\|^2\nonumber\\
&=f(x^{r+1}) +\frac{1}{2\beta^{r+1}}\left(\|\mu^{r+1}\|^2-\|\mu^r\|^2+ \|\mu^{r+1}-\mu^r\|^2\right)+\frac{\beta^{r+1}}{2}\|Ax^{r+1}-b\|^2\nonumber\\
&=f(x^{r+1}) +\frac{1}{2\beta^{r+1}}\|\mu^{r+1}\|^2-\frac{1}{2\beta^{r}}\|\mu^r\|^2+ \frac{1}{2\beta^{r+1}}\|\mu^{r+1}-\mu^r\|^2+ \left(\frac{1}{2\beta^{r}}-\frac{1}{2\beta^{r+1}}\right)\|\mu^r\|^2+\frac{\beta^{r+1}}{2}\|Ax^{r+1}-b\|^2\nonumber\\
&\ge f(x^{r+1}) +\frac{1}{2\beta^{r+1}}\|\mu^{r+1}\|^2-\frac{1}{2\beta^{r}}\|\mu^r\|^2+ \frac{1}{2\beta^{r+1}}\|\mu^{r+1}-\mu^r\|^2+\frac{\beta^{r+1}}{2}\|Ax^{r+1}-b\|^2\nonumber
\end{align}
where we have used the fact that $\beta^{r+1}\ge \beta^r$. Note that $t_0$ in \eqref{eq:clean} is a finite number hence $\frac{1}{2\beta^{t_0}}\|\mu^{t_0}\|^2$ is finite, and utilize Assumption [A2], we conclude that 
\begin{align}
\sum_{r={t_0}}^{\infty} &L_{\beta^{r+1}}(x^{r+1}, \mu^{r+1})>-\infty. 
\end{align}

Combining the above with the fact that the remaining terms of the potential function are all nonnegative, we conclude
\begin{align}\label{eq:P:bound:inc}
\sum_{r=1}^{\infty}P_{\beta^{r+1}, c}(x^{r+1}, x^r, \mu^{r+1}) >-\infty. 
\end{align}
Combining \eqref{eq:P:bound:inc} and the bound \eqref{eq:clean} (which is true for a finite $t_0>0$),  we conclude that the potential function $P_{\beta^{r+1}, c}(x^{r+1}, x^r, \mu^{r+1}) $ is lower bounded for all $r$.

\noindent{\bf Step 6.} In this step we show that the successive differences of various quantities converge. 

The lower boundedness of the potential function combined with the bound \eqref{eq:clean} (which is true for a finite $t_0>0$) implies that
\begin{subequations}
	\begin{align}
	&\sum_{r=1}^{\infty} \beta^{r+1} \|x^{r+1}-x^r\|^2<\infty. \label{eq:x:to:bound}\\
	&\sum_{r=1}^{\infty} (\beta^{r+1})^2 \left\|(x^{r+1}-x^r)-(x^{r}-x^{r-1})\right\|_{B^T B}^2<\infty\label{eq:x:diff:to:bound}.
	\end{align}
\end{subequations}
Therefore, we have
\begin{subequations}
	\begin{align}
	&\beta^{r+1} \|x^{r+1}-x^r\|^2\to 0. \label{eq:x:to:0}\\
	&\left (\beta^{r+1})^2\|(x^{r+1}-x^r)-(x^{r}-x^{r-1})\right\|_{B^T B}^2\to 0\label{eq:x:diff:to:0}.
	\end{align}
\end{subequations}
These two facts applied to \eqref{eq:differ:mu:inc}, combined with $\mu^{r+1}-\mu^r\in\mbox{col} (A)$, indicate that the following is true 
\begin{align}\label{eq:mu:converge:inc}
\mu^{r+1}-\mu^r\to 0.
\end{align}
Also \eqref{eq:clean} implies that the potential function is {\it upper bounded} as well, and this indicates that 
\begin{align}
\frac{c \beta^{r+1}\beta^r}{2}\|Ax^{r+1}-b\|^2 \; \mbox{is bounded}, \quad \frac{c \beta^{r+1}\beta^r}{2}\|x^r-x^{r+1}\|^2\; \mbox{is bounded}. 
\end{align} 

The second of the above inequality implies that $\beta^{r+1}B^T B(x^{r+1}-x^r)$ is bounded. If we further assume that $\nabla f(x)$ is bounded, and use \eqref{eq:mu:update:inc2}, we can conclude that $\{\mu^r\}$ is bounded.

\noindent{\bf Step 7.} Next we show that every limit point of $(x^r, \mu^r)$ converges to a stationary solution of problem \eqref{eq:original}. Let us pass a subsequence $\cK$ to $(x^r, \mu^r)$ and denote $(x^*, \mu^*)$ as its limit point. For notational simplicity, in the following the index $r$ all belongs to the set $\cK$. 

From relation \eqref{eq:x:to:bound} we have that any given $\epsilon>0$, there exists $t$ large enough the following is true
\begin{align}\label{eq:x:difference:epsilon}
\sum_{r=t-1}^{\infty}\beta^{r+1}\|x^{r+1}-x^r\|^2\le \frac{\epsilon}{c \omega 16}.
\end{align}

Utilizing \eqref{eq:def:v}, we have that the following is true
\begin{align}
\begin{split}\label{eq:v:summable}
&\sum_{r=1}^{\infty} \frac{1}{\beta^{r+1}} \|v^{r+1}\|^2<\infty, \quad  
\lim_{t\to\infty}\sum_{r=t}^{\infty}(\beta^{r+1})^2\left\|(x^{r+1}-x^r)-(x^{r}-x^{r-1})\right\|_{B^T B}^2= 0.
\end{split}
\end{align}
The first relation implies that $\lim\inf_{r\to\infty}\|v^{r+1}\|= 0$. Applying these relations to \eqref{eq:mu:bound:inc}, we have
\begin{align*}
\sum_{r=1}^{\infty}\frac{1}{\beta^{r+1}}\|\mu^{r+1}-\mu^r\|^2<\infty.
\end{align*}
This implies that for any given $\epsilon>0$, $c>0$, there exists an index $t$ sufficiently large such that 
\begin{align}\label{eq:mu:difference:bounded}
\sum_{r=t-1}^{\infty}\frac{1}{\beta^{r+1}}\|\mu^{r+1}-\mu^r\|^2<\frac{\epsilon^2}{4096  L \|B^T B\|\omega (1+c)}.
\end{align}
Applying this inequality and \eqref{eq:x:difference:epsilon} to \eqref{eq:clean}, we have that for large enough $t$ and for any $T>t$ the following is true 
\begin{align}\label{eq:one:step:descent}
P_{\beta^{T+1}, c}(x^{T+1}, x^T, \mu^{T+1}) - P_{\beta^{t}, c}(x^{t}, x^{t-1}, \mu^{t})\le -\sum_{r=t}^{T}\left(\frac{\beta^{r+1}}{48}\|x^{r+1}-x^r\|^2\right) + \frac{\epsilon^2}{4096L \|B^T B\|}.
\end{align}

Next we modify a classical argument in \cite[Proposition 3.5]{bertsekas96} to show that 
\begin{align*}
\lim_{r\to\infty}\|v^{r+1}\|\to 0.
\end{align*} 
We already know from the first relation in \eqref{eq:v:summable} that $\lim\inf_{r\to\infty}\|v^{r+1}\|= 0$. Suppose that $\|v^{r+1}\|$ does not converge to $0$, then we must have $\lim\sup_{r\to\infty} \|v^{r+1}\| >0$. Hence there exists an $\epsilon>0$ such that $\|v^{r+1}\|<\epsilon/2$ for infinitely many $r$, and $\|v^{r+1}\|>\epsilon$ for infinitely many $\epsilon$. Then there exists an infinite subset of iteration indices $\cR$ such that for each $r\in \cR$, there exits a $t(r)$ such that
\begin{align*}
&\|v^r\|<\epsilon/2, \quad, \|v^{t(r)}\|>\epsilon, \nonumber\\
& \epsilon/2<\|v^{t}\|\le \epsilon, \quad \forall~r<t<t(r).
\end{align*}
Using the fact that $\lim_{r\in \cK}\mu^{r} = \mu^*$, we have that for $r$ large enough, the following is true for all $t\ge 0$
\begin{align}\label{eq:key:mu}
\|\mu^{r} - \mu^{r+t}\|\le \frac{\epsilon}{8}\frac{1}{\|(B^T B)^{-1}\| \|A^T A\|}.
\end{align}
Without loss of generality we can assume that  this relation holds for all $r\in \cR$.
Note that the following is true
\begin{align}
\frac{\epsilon}{2}\le \|v^{t(r)}\|-\|v^r\|&\le \|v^{t(r)}-v^r\|= \left\|(B^T B)^{-1}\sum_{t=r}^{t(r)-1}\left(\nabla f(x^{t+1})-\nabla f(x^t)+A^T(\mu^{t+1}-\mu^t)\right)\right\|\nonumber\\
&\le \|(B^T B)^{-1}\|\left(\sum_{t=r}^{t(r)-1}\|\nabla f(x^{t+1})-\nabla f(x^t)\|+ \|A^T A\|\|\mu^{t(r)}-\mu^r\|\right)\nonumber\\
&\stackrel{\eqref{eq:def:v}}\le \|(B^T B)^{-1}\|\left(\sum_{t=r}^{t(r)-1}\frac{L}{\beta^{t+1}}\|v^{t+1}\|+ \|A^T A\|\|\mu^{t(r)}-\mu^r\|\right)\nonumber\\
&\le {\epsilon} L \|(B^T B)^{-1}\|\sum_{t=r}^{t(r)-1}\frac{1}{\beta^{t+1}}+\frac{\epsilon}{8}
\end{align}
where in the last inequality we have used \eqref{eq:key:mu} and the fact that for all $t\in(r+1, t(r))$, we have $\|v^{t}\|<{\epsilon}$. 
This implies that 
\begin{align}\label{eq:contradiction}
\frac{3}{8 L \|(B^T B)^{-1}\|} \le \sum_{t=r}^{t(r)-1}\frac{1}{\beta^{t+1}}.
\end{align}
Using the descent of the potential function \eqref{eq:one:step:descent} we have, for $r\in \cR$ and $r$ large enough
\begin{align}\label{eq:contradict}
&P_{\beta^{t(r)}, c}(x^{t(r)}, x^{t(r)-1}, \mu^{t(r)}) - P_{\beta^{r}, c}(x^{r}, x^{r-1}, \mu^{r})\nonumber\\
&\le-\sum_{t=r}^{t(r)-1}\frac{1}{48\beta^{t+1}}\|v^{t+1}\|^2 +\frac{\epsilon^2}{4096 L \|B^T B\|}\nonumber\\
&\stackrel{\rm (i)}\le -\left(\frac{\epsilon}{4}\right)^2\sum_{t=r}^{t(r)-1}\frac{1}{48\beta^{t+1}}+\frac{\epsilon^2}{4096 L \|B^T B\|}\nonumber\\
&\stackrel{\rm (ii)}\le -\frac{\epsilon^2}{2048L\|B^T B\|}+\frac{\epsilon^2}{4096L \|B^T B\|}\nonumber\\
&\le -\frac{\epsilon^2}{4096 L \|B^T B\|}
\end{align}
where in ${\rm (i)}$ we have used the fact that for all $r\in \cR$, $\|v^{r+i}\|\ge\frac{\epsilon}{2}$ for $i=1,\cdots, t(r)$; in ${\rm (ii)}$ we have used \eqref{eq:contradiction}.
However we know that the potential function is converging, i.e.,
$$\lim_{r\to\infty}P_{\beta^{t(r)}, c}(x^{t(r)}, x^{t(r)-1}, \mu^{t(r)}) \to P_{\beta^{r}, c}(x^{r}, x^{r-1}, \mu^{r})=0$$ 
which contradicts to \eqref{eq:contradict}. Therefore we conclude that $\|v^{r+1}\|\to 0$.  

Finally, combining $\|v^{r+1}\|\to 0$ with the convergence of $\mu^{r+1}-\mu^r$ (cf. \eqref{eq:mu:converge:inc}), we conclude that every limit point of $\{x^r, \mu^r\}$ satisfies
\begin{align*}
\nabla f(x^*)+ A^T \mu^* =0, \quad A x^* =b.
\end{align*}
Therefore it is a stationary solution for problem \eqref{eq:original}. This completes the proof.

\section{Proof of Convergence for Algorithm 6} \label{sec:proof:matrix}

To make the derivation compact, define the following matrix 
\begin{align}\label{eq:M}
\bfM^{r+1} &:=\nabla_{\bX}f(\bX^{r+1}, Y^{r+1}) \nonumber\\
&=\left[((X^{r+1}_1y^{r+1}_1)-z_1)(y^{r+1}_1)^T+ 2\gamma X^{r+1}_1; \cdots; ((X^{r+1}_N y^{r+1}_N)-z_N)(y^{r+1}_N)^T+ 2\gamma X^{r+1}_N\right] .
\end{align}
The proof consists of six steps. 

\noindent{\bf Step 1}. First we note that the optimality condition for the $\bX$-subproblem \eqref{eq:x:update:matrix} is given by
\begin{align}\label{eq:omega:expression}
\bA^T \bfOmega^{r+1} = -\bfM^{r+1} -\beta \langle \bB^T \bB (\bX^{r+1}-\bX^{r})\rangle 
\end{align}
By utilizing the fact that that $\bfOmega^{r+1}-\bfOmega^r$ lies in the column space of $\bA$, and the eigenvalues of $A^TA$ equal to the eigenvalue of $\bA^T \bA$, we have the following bound 
\begin{align*}
\|\bfOmega^{r+1}-\bfOmega^r\|_F^2\le \frac{2}{{\sigma_{\min}(A^T A)}}\left(\|\bfM^{r+1}-\bfM^r\|_F^2 +\beta^2\|\bB^T \bB[(\bX^{r+1}-\bX^{r}) - (\bX^{r}-\bX^{r-1})]\|_F^2\right).
\end{align*}
Next  let us analyze the first term in the rhs of the above inequality. The following identity holds true
\begin{align}\label{eq:Mi}
&\|\bfM^{r+1}-\bfM^r\|_F^2\nonumber\\
&=\sum_{i=1}^{N}\|(X^{r+1}_i y^{r+1}_i - z_i)(y^{r+1}_i)^T-(X^{r}_i y^r_i - z_i)(y^{r}_i)^T + 2\gamma (X^{r+1}_i-X^{r}_i)\|^2_F\nonumber\\
&\le \sum_{i=1}^{N}4\|X^{r+1}_i-X^{r}_i\|^2_F \|y^{r+1}_i (y^{r+1}_i)^T\|^2+4\|X^r_i y^r_i - z_i\|^2\|y^{r+1}_i-y^r_i\|^2 +4\|X^r_i(y^{r+1}_i-y^r_i)\|^2\|y^{r+1}_i\|^2\nonumber\\
&\quad + 16\gamma^2\|X^{r+1}_i-X^{r}_i\|^2_F \nonumber\\
&\le  \sum_{i=1}^{N}4(\tau^2+4\gamma^2)\|X^{r+1}_i-X^r_i\|_F^2+4{\theta^r_i}\|y^{r+1}_i-y^r_i\|^2 + 4{\tau}\|X^r_i(y^{r+1}_i-y^r_i)\|^2
\end{align}
where in the last inequality we have defined the constant $\theta^r_i$ as 
\begin{align}\label{eq:theta}
\theta^r_i:= \|X^r_i y^r_i - z_i\|^2.
\end{align}

Therefore, combining the above two inequalities, we obtain
\begin{align}\label{eq:omega:difference}
&\frac{1}{\beta}\|\bfOmega^{r+1}-\bfOmega^r\|_F^2\nonumber\\
&\le \frac{8}{\beta{\sigma_{\min}(A^T A)}}\sum_{i=1}^{N} \left((\tau^2+4\gamma^2) \|X^{r+1}_i-X^r_i\|_F^2 +{\theta^r_i}\|y^{r+1}_i-y^r_i\|^2 + {\tau}\|X^r_i(y^{r+1}_i-y^r_i)\|^2\right) \nonumber\\
&\quad +\frac{2\beta}{{\sigma_{\min}(A^T A)}}\|\bB^T \bB[(\bX^{r+1}-\bX^{r}) - (\bX^{r}-\bX^{r-1})]\|_F^2
\end{align}

\noindent{\bf Step 2}. Next let us analyze the descent of the augmented Lagrangian. First we have
\begin{align}\label{eq:y:descent:matrix}
&L_{\beta}(\bX^{r}, Y^{r+1},\bfOmega^{r})-L_{\beta}(\bX^{r}, Y^{r},\bfOmega^{r})\nonumber\\
&=\sum_{i=1}^{N}\left(\frac{1}{2}\|X^r_i y^{r+1}_i-z_i\|^2+h_i(y^{r+1}_i) - \frac{1}{2}\|X^r_i y^{r}_i-z_i\|^2-h_i(y^r_i)\right)\nonumber\\
&\le \sum_{i=1}^{N}\left(\frac{1}{2}\|X^r_i y^{r+1}_i-z_i\|^2+h_i(y^{r+1}_i)+\frac{\theta^r_i}{2}\|y_i^{r+1}-y_i^r\|^2 - \frac{1}{2}\|X^r_i y^{r}_i-z_i\|^2-h_i(y^r_i)\right)\nonumber\\
& \le  \sum_{i=1}^{N}\bigg(\left\langle (X^r_i)^T(X^r_i y^{r+1}_i-z_i)+\theta^r_i(y^{r+1}_i-y^r_i), y^{r+1}_i-y^r_i\right\rangle -\frac{1}{2}\|X^r_i(y^{r+1}_i-y^r_i)\|^2-\frac{\theta^r_i}{2}\|y_i^{r+1}-y_i^r\|^2 \nonumber\\
&\quad +\langle \zeta^{r+1}_i, y^{r+1}_i-y^r_i\rangle \bigg)\nonumber\\
&\le - \sum_{i=1}^{N}\left( \frac{1}{2}\|X^r_i(y^{r+1}_i-y^r_i)\|^2+\frac{\theta^r_i}{2}\|y_i^{r+1}-y_i^r\|^2 \right)
\end{align}
where in the second to the last equality we have used the convexity of $h_i$, and $\zeta^{r+1}_i\in \partial h_i(y_i^{r+1})$; the last inequality uses the optimality condition of the $y$-step \eqref{eq:y:update:matrix}. Similarly, we can show that 
\begin{align}\label{eq:x:descent:matrix}
&L_{\beta}(\bX^{r+1}, Y^{r+1},\bfOmega^{r})-L_{\beta}(\bX^{r+1}, Y^{r},\bfOmega^{r})\le - \frac{\beta+2\gamma}{2}\|\bX^{r+1}-\bX^r\|^2_F
\end{align}
where we have utilized the fact that $\bA^T\bA+\bB^T \bB =  2\mathbf{D}\succeq {\bf I}_{NM}.$ Therefore, combining the estimate \eqref{eq:omega:difference}, we obtain
\begin{align}\label{eq:lag:descent:matrix}
&L_{\beta}(\bX^{r+1}, Y^{r+1},\bfOmega^{r+1})-L_{\beta}(\bX^{r}, Y^{r},\bfOmega^{r})\nonumber\\
&\le -\left(\frac{\beta+2\gamma}{2}-\frac{8(\tau^2+4\gamma^2)}{\beta{\sigma_{\min}(A^T A)}}\right)\sum_{i=1}^{N} \|X^{r+1}_i-X^r_i\|_F^2 -\sum_{i=1}^{N}\left(\frac{\theta^r_i}{2}-\frac{8\theta^r_i}{\beta\sigma_{\min}(A^T A)}\right)\|y^{r+1}_i-y^r_i\|^2 \nonumber\\
&\quad -\left(\frac{1}{2}-\frac{{8\tau}}{\sigma_{\min}(A^T A)\beta}\right)\sum_{i=1}^{N}\|X^r_i(y^{r+1}_i-y^r_i)\|^2 \nonumber\\
&\quad +\frac{2\beta}{{\sigma_{\min}(A^T A)}}\|\bB^T \bB[(\bX^{r+1}-\bX^{r}) - (\bX^{r}-\bX^{r-1})]\|_F^2.
\end{align}

\noindent {\bf Step 3. } This step follows  Lemma \ref{lemma:second:descent} in the analysis of Algorithm 1. In particular, after writing down the optimality condition of the $\bX^{r+1}$ and $\bX^{r}$ step, we can obtain
\begin{align*}
&\langle \bA^T (\bfOmega^{r+1}-\bfOmega^r), \bX^{r+1}-\bX^r\rangle \nonumber\\
&\le -\left\langle  \bfM^{r+1}-\bfM^r +\beta\bB^{T}\bB \left[(\bX^{r+1}-\bX^r) - (\bX^{r}-\bX^{r-1})\right]   , \bX^{r+1}-\bX^r\right\rangle .
\end{align*}

Then it is easy to show that the above inequality implies the following (which utilizes the convexity of  $h$)
\begin{align}
&\frac{\beta}{2}\left(\langle \bA\bX^{r+1}, \bA\bX^{r+1}\rangle +\langle\bB^T\bB(\bX^{r+1}-\bX^r), \bX^{r+1}-\bX^r \rangle \right)\nonumber\\
&\le  \frac{\beta}{2}\left(\langle \bA\bX^{r}, \bA\bX^{r}\rangle+\langle\bB^T\bB(\bX^{r}-\bX^{r-1}), \bX^{r}-\bX^{r-1} \rangle\right) -\frac{\beta}{2}\langle \bA(\bX^{r+1}-\bX^r), \bA(\bX^{r+1}-\bX^r)\rangle  \nonumber\\
&\quad -\langle \bfM^{r+1}-\bfM^r, \bX^{r+1}-\bX^r\rangle -\frac{\beta}{2}\|\bB[(\bX^{r+1}-\bX^{r}) - (\bX^{r}-\bX^{r-1})]\|_F^2 .\nonumber
\end{align}

Note the following fact
\begin{align}\label{eq:diff:M}
&-\langle \bfM^{r+1}-\bfM^r, \bX^{r+1}-\bX^r\rangle \nonumber\\
&= -\langle \nabla_\bX f(\bX^{r+1}, y^{r+1}) - \nabla_\bX f(\bX^{r}, y^{r}) , \bX^{r+1}-\bX^r\rangle\nonumber\\
&= -\langle \nabla_\bX f(\bX^{r+1}, y^{r+1})- \nabla_\bX f(\bX^{r}, y^{r+1}) + \nabla_\bX f(\bX^{r}, y^{r+1})  - \nabla_\bX f(\bX^{r}, y^{r}) , \bX^{r+1}-\bX^r\rangle\nonumber\\
&\stackrel{\rm(i)}\le -\langle \nabla_\bX f(\bX^{r}, y^{r+1})  - \nabla_\bX f(\bX^{r}, y^{r}) , \bX^{r+1}-\bX^r\rangle \nonumber\\
&\stackrel{\rm(ii)}\le \frac{1}{2d}\|\nabla_\bX f(\bX^{r}, y^{r+1})  - \nabla_\bX f(\bX^{r}, y^{r})\|_F^2+\frac{d}{2}\|\bX^{r+1}-\bX^r\|_F^2\nonumber\\
&\stackrel{\rm(iii)}\le \frac{1}{d}\sum_{i=1}^{N}\left({\theta^r_i}\|y^{r+1}_i-y^r_i\|^2 + {\tau}\|X^r_i(y^{r+1}_i-y^r_i)\|^2\right)+\frac{d}{2}\|\bX^{r+1}-\bX^r\|_F^2
\end{align}
where in ${\rm (i)}$ we utilize the convexity of $f(\bX, y)$ wrt $\bX$ for any fixed $y$; in ${\rm(ii)}$ we use the Cauchy-Swarch inequality, where $d>0$ is a constant (to be determined later); ${\rm(iii)}$ is true due to a similar calculation as in \eqref{eq:Mi}.
	 
Therefore overall we have
\begin{align}\label{eq:descent:matrix:2}
&\frac{\beta}{2}\left(\langle \bA\bX^{r+1}, \bA\bX^{r+1}\rangle +\langle\bB^T\bB(\bX^{r+1}-\bX^r), \bX^{r+1}-\bX^r \rangle \right)\nonumber\\
&\le  \frac{\beta}{2}\left(\langle \bA\bX^{r}, \bA\bX^{r}\rangle+\langle\bB^T\bB(\bX^{r}-\bX^{r-1}), \bX^{r}-\bX^{r-1} \rangle\right) -\frac{\beta}{2}\langle \bA(\bX^{r+1}-\bX^r), \bA(\bX^{r+1}-\bX^r)\rangle \nonumber\\
&\quad +\frac{1}{d}\sum_{i=1}^{N}\left({\theta^r_i}\|y^{r+1}_i-y^r_i\|^2 + {\tau}\|X^r_i(y^{r+1}_i-y^r_i)\|^2\right)+\frac{d}{2} \|\bX^{r+1}-\bX^r\|^2_F\nonumber\\
&\quad -\frac{\beta}{2}\|\bB[(\bX^{r+1}-\bX^{r}) - (\bX^{r}-\bX^{r-1})]\|_F^2 
\end{align}
\noindent {\bf Step 4. } Let us define the potential function as 
\begin{align}\label{eq:potential:matrix}
&P_{\beta, c}(\bX^{r+1}, \bX^{r}, Y^{r+1}, \bfOmega^{r+1})\nonumber\\
&:=L_{\beta}(\bX^{r+1}, Y^{r+1},\bfOmega^{r+1}) + \frac{c\beta}{2}\left(\langle \bA\bX^{r+1}, \bA\bX^{r+1}\rangle +\langle\bB^T\bB(\bX^{r+1}-\bX^r), \bX^{r+1}-\bX^r \rangle \right).
\end{align}
Then utilize the bound in \eqref{eq:lag:descent:matrix} in Step 2 and bounds \eqref{eq:descent:matrix:2}  in Step 3, we obtain
\begin{align*}
&P_{\beta, c}(\bX^{r+1}, \bX^{r}, Y^{r+1}, \bfOmega^{r+1})-P_{\beta, c}(\bX^{r}, \bX^{r-1}, Y^{r}, \bfOmega^{r})\nonumber\\
&\le\quad -\left(\frac{\beta+2\gamma}{2}-\frac{8(\tau^2+4\gamma^2)}{\beta{\sigma_{\min}(A^T A)}}-\frac{cd}{2}\right)\sum_{i=1}^{N} \|X^{r+1}_i-X^r_i\|_F^2\nonumber\\ &-\sum_{i=1}^{N}\left(\frac{\theta^r_i}{2}-\frac{8\theta^r_i}{\beta\sigma_{\min}(A^T A)}-\frac{c\theta^r_i}{d}\right)\|y^{r+1}_i-y^r_i\|^2 \nonumber\\
&\quad -\left(\frac{1}{2}-\frac{{8\tau}}{\sigma_{\min}(A^T A)\beta}-\frac{c\tau}{d}\right)\sum_{i=1}^{N}\|X^r_i(y^{r+1}_i-y^r_i)\|^2 \nonumber\\
&\quad -\left(\frac{c\beta}{2}-\frac{2\beta\|B^T B\|}{{\sigma_{\min}(A^T A)}}\right)\|\bB[(\bX^{r+1}-\bX^{r}) - (\bX^{r}-\bX^{r-1})]\|_F^2.
\end{align*}

Therefore the following are the condition that guarantees the descent of the potential function
\begin{align}\label{eq:condition:matrix}
\begin{split}
\frac{\beta+2\gamma}{2}-\frac{8(\tau^2+4\gamma^2)}{\beta{\sigma_{\min}(A^T A)}}-\frac{cd}{2}&>0, \quad \frac{1}{2}-\frac{{8}}{\sigma_{\min}(A^T A)\beta}-\frac{c}{d}>0\\
\frac{1}{2}-\frac{{8\tau}}{\sigma_{\min}(A^T A)\beta}-\frac{c\tau}{d}&>0, \quad  \frac{c\beta}{2}-\frac{2\beta\|B^T B\|}{{\sigma_{\min}(A^T A)}} >0.
\end{split}
\end{align}
To see that it is always possible to find the tuple $(\beta, c, d)$, first let us set $c$  such that the last inequality is satisfied
\begin{align}
c>\frac{4\|B^T B\|}{\sigma_{\min}(A^T A)}\label{eq:c:bound:matrix}.
\end{align}
Second, let us pick any $d$ such that the following is true
\begin{align*}
d>\max\{ 2 c\tau, 2c\}.
\end{align*}
Then clearly it is possible to make $\beta$ large enough such that all the four conditions in \eqref{eq:condition:matrix} are satisfied.

\noindent {\bf Step 5. }  We need to prove that the potential function is lower bounded. We lower bound the  augmented Lagrangian as follows 
\begin{align}\label{eq:lag:matrix}
&L_{\beta}(\bX^{r+1}, Y^{r+1}, \bfOmega^{r+1}) \nonumber\\
&= \sum_{i=1}^{N}\frac{1}{2}\|X^{r+1}_i y^{r+1}_i - z_i\|^2+\gamma\|X^{r+1}_i\|^2_F+h_i(y^{r+1}_i)+\langle \bfOmega^{r+1}, \bA\bX^{r+1}\rangle  +\frac{\beta}{2}\langle \bA\bX^{r+1}, \bA\bX^{r+1}\rangle \nonumber\\
&= \sum_{i=1}^{N}\frac{1}{2}\|X^{r+1}_i y^{r+1}_i - z_i\|^2+\gamma\|X^{r+1}_i\|^2_F+h_i(y^{r+1}_i)  +\frac{\beta}{2}\langle \bA\bX^{r+1}, \bA\bX^{r+1}\rangle \nonumber\\
&\quad +\frac{1}{2\beta}\left(\|\bfOmega^{r+1}-\bfOmega^r\|_F^2+\|\bfOmega^{r+1}\|_F^2-\|\bfOmega^{r}\|_F^2\right).
\end{align}
Then by the same argument leading to \eqref{eq:lower:bound:case1}, we conclude that as long as $h_i$ is lower bounded over its domain, then the potential function will be lower bounded. 

{\bf Step 6.} Combining the result in Step 5 and Step 4, we conclude the following
\begin{subequations}\label{eq:converge}
	\begin{align}
	&\sum_{i=1}^{N}\|X^{r+1}_i-X^r_i\|_F^2\to 0, \quad \sum_{i=1}^{N}\|y_i^{r+1}-y^r_i\|^2\to 0\\
	&\sum_{i=1}^{N}\|X^r_i(y_i^{r+1}-y^r_i)\|^2\to 0, \quad \left\|\bB^T \bB [(\bX^{r+1}-\bX^r)-(\bX^r-\bX^{r-1})]\right\|_F\to 0.
	\end{align}
\end{subequations}

Then utilizing \eqref{eq:omega:difference}, we have 
\begin{align*}
\bfOmega^{r+1}-\bfOmega^r \to\mathbf{0},\; \mbox{or equivalently}\; \bA\bX^{r+1} \to 0.
\end{align*}
That is, in the limit the network-wide consensus is achieved. Next we show that the primal and dual iterates are bounded. 

Note that the potential function is both lower and upper bounded. Combined with \eqref{eq:converge} we must have that the augmented Lagrangian is both upper and lower bounded. Using the expression \eqref{eq:lag:matrix},  the assumption that $h_i (y_i)$ is lower bounded, and the fact that $y_i$ is bounded, we have that in the limit, the following term is bounded
$$\sum_{i=1}^{N}\frac{1}{2}\|X^{r+1}_i y^{r+1}_i - z_i\|^2+\gamma\|X^{r+1}_i\|_F^2.$$
This implies that the primal variable sequence $\{X^{r+1}_i\}$ are bounded for all $i$. To show the boundedness of the dual sequence, note that $\bfOmega^{r+1}\in\mbox{col}(\bA)$ (due to the initialization that $\bfOmega^0=\mathbf{0}$). Therefore using \eqref{eq:omega:expression} we have
\begin{align*}
\sigma_{\min}(\bA^T \bA)\|\bfOmega^{r+1}\|_F^2 \le  2\|\bfM^{r+1}\|_F^2 +2 \beta\|  \bB^T \bB (\bX^{r+1}-\bX^{r}) \|_F^2
\end{align*}
Note that from the expression of $\bfM$ in \eqref{eq:M}, we see that $\{\bfM^{r+1}\}$ is bounded because both $\bX^{r+1}$ and $Y^{r+1}$ are bounded. Similarly, the second term on the rhs of the above inequality is bounded because $\bX^{r+1}\to \bX^r$. These two facts imply that $\{\bfOmega^{r+1}\}$ is bounded as well. 

Arguing the convergence to stationary point as well as the convergence rate follows exact the same steps as in the proof of Theorem \ref{thm:main}.

\newpage

\bibliographystyle{IEEEbib}

\bibliography{ref,biblio,distributed_opt}

\end{document}